\documentclass[12pt]{article}
\usepackage{amssymb}
\usepackage{rotating}
\usepackage{latexsym}

\newtheorem{conj}{Conjecture}[section]
\newtheorem{cor}{Corollary}[section]
\newtheorem{defn}{Definition}[section]
\newtheorem{hyp}{Hypothesis}[section]
\newtheorem{lemma}{Lemma}[section]
\newtheorem{prop}{Proposition}[section]

\newtheorem{thm}{Theorem}[section]

\newcounter{example}[section]
\newcounter{rem}[section]

\newcommand{\bbC}{\ensuremath{\mathbb{C}}}
\newcommand{\bbCp}{\ensuremath{\mathbb{C}_p}}
\newcommand{\bbN}{\ensuremath{\mathbb{N}}}

\newcommand{\bbQ}{\ensuremath{\mathbb{Q}}}

\newcommand{\bbR}{\ensuremath{\mathbb{R}}}
\newcommand{\bbZ}{\ensuremath{\mathbb{Z}}}
\newcommand{\bbZp}{\ensuremath{\mathbb{Z}_p}}

\newcommand{\beql}[1]{\begin{equation}\label{#1}}
\newcommand{\bPf}{\noindent \textsc{Proof\ }}

\newcommand{\bt}{\ensuremath{{\bf t}}}

\newcommand{\cA}{\ensuremath{\mathcal{A}}}
\newcommand{\cB}{\ensuremath{\mathcal{B}}}
\newcommand{\cD}{\ensuremath{\mathcal{D}}}
\newcommand{\cE}{\ensuremath{\mathcal{E}}}
\newcommand{\cF}{\ensuremath{\mathcal{F}}}
\newcommand{\cf}{\emph{cf.}}
\newcommand{\cI}{\ensuremath{\mathcal{I}}}

\newcommand{\Cl}{{\rm Cl}}

\newcommand{\cM}{\ensuremath{\mathcal{M}}}
\newcommand{\cN}{\ensuremath{\mathcal{N}}}
\newcommand{\cO}{\ensuremath{\mathcal{O}}}

\newcommand{\cR}{\ensuremath{\mathcal{R}}}
\newcommand{\cS}{\ensuremath{\mathcal{S}}}

\newcommand{\displaymapdef}[5]
{\[
\begin{array}{rcrcl}
 #1 &:& #2 &\longrightarrow& #3 \\
    & &    &                    \\
    & & #4 &\longmapsto    & #5
\end{array}
\]}

\newcommand{\eeq}{\end{equation}}
\newcommand{\eg}{\emph{e.g.}}
\newcommand{\ePf}{\hspace*{\fill}~$\Box$\vertsp\par}
\newcommand{\etc}{\emph{etc.}}
\newcommand{\fa}{\ensuremath{\mathfrak{a}}}

\newcommand{\fc}{\ensuremath{\mathfrak{c}}}
\newcommand{\ff}{\ensuremath{\mathfrak{f}}}

\newcommand{\fm}{\ensuremath{\mathfrak{m}}}

\newcommand{\fp}{\ensuremath{\mathfrak{p}}}

\newcommand{\fq}{\ensuremath{\mathfrak{q}}}

\newcommand{\fw}{\ensuremath{\mathfrak{w}}}
\newcommand{\fW}{\ensuremath{\mathfrak{W}}}
\newcommand{\fz}{\ensuremath{\mathfrak{z}}}
\newcommand{\Gal}{{\rm Gal}}

\newcommand{\ie}{\emph{i.e.}}
\newcommand{\inv}{^{-1}}
\newcommand{\ndiv}{\nmid}
\newcommand{\nin}{\not\in}
\newcommand{\ord}{{\rm ord}}

\newcommand{\refeq}[1]{~(\ref{eq:#1})}
\newcommand{\rem}{\refstepcounter{rem}\noindent{\sc Remark \therem}}

\newcommand{\sgn}{{\rm sgn}}
\renewcommand{\theexample}{\thesection.\arabic{example}}
\renewcommand{\therem}{\thesection.\arabic{rem}}

\newcommand{\vertsp}{\vspace{1ex}}

\newcommand{\barbbQ}{\ensuremath{\bar{\bbQ}}}
\newcommand{\Clfk}{\ensuremath{\Cl_\ff(k)}}

\newcommand{\Clmk}{\ensuremath{\Cl_\fm(k)}}

\newcommand{\ContZps}{{\rm Cont}(\bbZp^2,\bbC_p)}

\newcommand{\GQbarQ}{\ensuremath{\Gal(\barbbQ/\bbQ)}}

\newcommand{\kps}{\ensuremath{k_+^\times}}
\newcommand{\kzs}{\ensuremath{k_\fz^\times}}
\newcommand{\MeasZps}{{\rm Meas}(\bbZp^2,\bbC_p)}

\newcommand{\PhifT}{\ensuremath{\Phi_{\ff,T}}}

\newcommand{\PhifTpp}{\ensuremath{\Phi_{\ff,T_p,p}}}

\newcommand{\PhimT}{\ensuremath{\Phi_{\fm,T}}}
\newcommand{\PhimTp}{\ensuremath{\Phi_{\fm,T,p}}}

\newcommand{\QG}{\ensuremath{\bbQ G}}

\newcommand{\simz}{\sim_{\raise -3pt\hbox{$\scriptstyle
\kern -6pt \fz \kern 3pt$}}}
\newcommand{\Sr}{\ensuremath{^{[S,r]}}}
\newcommand{\St}{\ensuremath{^{[S,2]}}}
\newcommand{\tchi}{\ensuremath{{\tilde{\chi}}}}

\newcommand{\teSgt}{\ensuremath{\tilde{e}_{S,>2}}}

\newcommand{\tiota}{\ensuremath{{\tilde{\iota}}}}

\newcommand{\TrkQ}{{\rm Tr}_{k/\bbQ}}
\newcommand{\ta}{\tilde{a}}
\newcommand{\tA}{\tilde{A}}
\newcommand{\trho}{\tilde{\rho}}

\newcommand{\uiota}{\ensuremath{\underline{\iota}}}
\newcommand{\uX}{\ensuremath{\underline{X}}}

\newcommand{\WedrQUS}{\ensuremath{\textstyle\bigwedge_{\bbQ G}^r\bbQ U_S}}

\newcommand{\wmz}{\ensuremath{\fw_\fm^0}}

\newcommand{\WedtQUS}{\ensuremath{\textstyle\bigwedge_{\bbQ G}^2\bbQ U_S}}
\newcommand{\WedtZUS}{\ensuremath{\textstyle\bigwedge_{\bbZ G}^2 U_S}}
\newcommand{\ZG}{\ensuremath{\bbZ G}}

\newcommand{\ZT}[2]{\ensuremath{Z_T(#1;#2)}}

\newcommand{\ZTpj}[2]{\ensuremath{Z_{T,p}^{(j)}(#1;#2)}}

\newcommand{\onecol}[2]{\multicolumn{#1}{|l|}{\qquad #2}}
\newcommand{\onecolc}[2]{\multicolumn{#1}{|c|}{#2}}
\newcommand{\onecoll}[2]{\multicolumn{#1}{|l|}{#2}}
\newcommand{\onecolr}[2]{\multicolumn{#1}{|r|}{#2}}

\newcommand{\phv}[2]{\vrule height #1pt depth #2pt width 0pt}
\begin{document}

%
%
%
%
%
%
%

\title{Verifying a $p$-Adic Abelian Stark Conjecture
at $s = 1$}
\author{X.-F. Roblot \\
IGD, Universit\'e Lyon I
\and
D. Solomon
\thanks{supported by an Advanced Fellowship from the EPSRC.}\\
King's College London}
\maketitle

\begin{abstract}
In a previous paper \cite{z1}, the second author developed a new
approach to the abelian $p$-adic Stark conjecture at $s=1$ and stated
related conjectures. The aim of the present paper is to develop and
apply techniques to numerically investigate one of these -- the `Weak
Refined Combined Conjecture' -- in fifteen cases.
\end{abstract}

\section{Introduction}
In the 1970's and 80's Harold Stark~\cite{St} made a series of
conjectures concerning the values at $s=1$ and $s=0$ of complex
Artin $L$-series attached to a Galois extension of number fields
$K/k$.  Subsequently, much theoretical and computational work has
been done, extending and testing these conjectures, with
particular attention paid recently to certain refined conjectures
in the case where $K/k$ is abelian (\cite{Ru}, \cite{Pop}).

In \cite{z1}, a new approach to the abelian case of the $p$-adic
conjecture at $s=1$ was developed and several related conjectures
were stated. The main aim of the present paper is to 
develop and apply techniques to numerically investigate one of
these
 -- the `Weak Refined Combined Conjecture'
(Conjecture~3.6 of~\cite{z1}, here
Conjecture~\ref{conj:2F}) -- in a
number of cases.

In Section~\ref{sec:conj}, we shall recall the definitions of the
complex and $p$-adic `twisted zeta functions'. They depend on two
parameters: a proper ideal $\ff$ of $\cO_k$, and a set $T$ of primes
ideals of $\cO_k$ (which, for the purpose of $p$-adic interpolation,
must contain the primes above $p$). Then the statements of the two
`combined conjectures' of~\cite{z1} are given. (The term `combined'
refers to the fact that each conjecture predicts both a complex and a
$p$-adic equality.) The main reference for this section is, of course,
\cite{z1}, but also \cite{twizas} which contains a reformulation
developed by the second author of a refined complex abelian conjecture
at $s=0$ originally made by Rubin in~\cite{Ru}. Briefly, the `Weak
Refined Combined Conjecture' takes the following form: we assemble all
the complex (\emph{resp.}\ $p$-adic) twisted zeta-functions for given
$\ff$ and $T$ into a single group-ring-valued function
$\Phi_{\ff,T}(s)$ (\emph{resp.}\ $\Phi_{\ff,T,p}(s)$, assuming that
$T$ contains the primes above $p$). Then, assuming that the primes in
$T$ do not divide $\ff$, the value of the latter at $s=1$ is
conjectured to be equal to the complex (\emph{resp.}\ $p$-adic)
group-ring-valued regulator of a certain element $\eta_{\ff,T}$
multiplied by an explicit algebraic constant. 
The element $\eta_{\ff,T}$ is constructed from certain $S$-units of 
the field $K$ which in this case is simply the ray-class field $k(\ff)$. 

Section~\ref{quadformula} develops a new formula 
to compute the element $\Phi_{\ff,T,p}(1)$.  We
concentrate on the case where $k$ is real quadratic although our
technique should extend to other totally real fields.  Relying as it does on
Shintani's method and the theory of $p$-adic measures, this technique
is very different in nature from that used to evaluate complex
$L$-functions. (For the latter we use~\cite{DT}.)  We stress that it
passes most naturally not by the analogous $p$-adic $L$-functions but
by the $p$-adic twisted zeta-functions themselves. Indeed, this was
one of the major reasons for introducing these functions and, in
preparation, their complex analogues.

Finally, Section~\ref{vericonj} is devoted to the numerical
investigation of the `Weak Refined Combined Conjecture' over a real
quadratic field. We first explain some procedures (for
example a continued fraction method based on ideas of Zagier) 
that greatly shorten the
calculation of $\Phi_{\ff,T,p}(1)$ using the formula of the
previous section. Then we explain the basis of our method for
verifying the conjecture. Since $[k:\bbQ]=2$ and $K$ is totally real,
our conjectures are `second order' in the sense that the relevant
complex $L$-functions have at least a double zero at $s=0$. The
corresponding fact on the `other side' of the conjectures is that both
the complex and $p$-adic regulators must be of rank $2$. One
consequence is that, unlike verifications of the (complex) first order
abelian Stark Conjectures (see for example~\cite{Ro}), the regulators
themselves do not determine $S$-units of $K$.
We therefore need different methods
for finding $K$ and $\eta_{\ff,T}$
and new criteria for affirming that
the latter satisfies the combined conjecture to the precision of our
computations.  
In fact, we use the methods of~\cite{Ro} (which actually
rely on the first-order complex conjecture!) to independently and verifiably
construct the ray-class field $K$.  We then illustrate our methodology
by numerically confirming the conjecture in fifteen different
examples, using a number of different primes $p$ in each example.  The
resulting data are displayed in tables at the end of the paper.  We
hope that they will serve to stimulate further interest in these
conjectures, their possible refinements and extensions.
\section{The $p$-adic Stark conjectures at $s=1$}\label{sec:conj}
The main reference for this section are \cite{twizas} (for the complex
twisted functions) and \cite{z1}.
\subsection{Complex twisted zeta functions}
Let $k\subset\bar{\bbQ}\subset\bbC$ be any number field of finite
degree over $\bbQ$ and let $\cO$ its ring of integers. Let $I$ be any
fractional ideal of $k$ and $\xi$ any character on (the additive group
of) $I$ with values in $\mu(\bbC)$, the complex roots of unity.  The
\emph{annihilator} of $\xi$ is the ideal $\ff\lhd\cO$ given by
$\ff=\{a\in\cO\,:\,\xi(ax)=1\ \forall\,x\in I\}$. Suppose that $\fz$
is the formal product of some subset of the real places of $k$ and
write $\fm$ for the \emph{cycle} that is the formal product $\ff\fz$.
We denote by $E_\fm$ the subgroup of finite index in $E(K):=\cO^\times$
consisting of the units that are congruent to $1$ modulo $\fm$ in the
usual sense.  For any finite set $T$ of finite places (prime ideals)
of $\cO$, the group $E_\fm$ acts by multiplication on the following
subset of $I$
\[
\cS(I,\fz,T):=\{a\in I\,:\,\mbox{$a\in\kzs$ and $(aI\inv,T)=1$}\}
\]
where $\kzs$ denotes the elements of $k^\times$ which are positive at
all places dividing $\fz$ and the notation $(J,T)=1$ indicates that an
ideal $J$ of $\cO$ has support disjoint from $T$. For $s\in\bbC,\
\Re(s)>1$ we consider the absolutely convergent Dirichlet series,
called the `twisted zeta-function' for these data, defined by
\begin{eqnarray}
Z_T(s;\xi,I,\fm):=\sum_{a\in\cS(I,\fz,T)/E_\fm}\frac{\xi(a)}{|I:(a)|^s}
               &=&
\sum_{a\in\cS(I,\fz,T)/E_\fm}\frac{\xi(a)}{N(aI\inv)^s}\nonumber\\
               &=&
NI^s\sum_{a\in\cS(I,\fz,T)/E_\fm}\xi(a)|{\rm N}_{k/\bbQ}(a)|^{-s}
\label{eq:1A}
\end{eqnarray}

Let $W_\ff$ be the set of all pairs
$(\psi,J)$, where $\psi$ is a character of annihilator $\ff$ on
a fractional ideal $J$. In~\cite{z1} a natural equivalence
relation (depending on $\fz$) was defined on $W_\ff$ in such a way that
$Z_T(s;\xi,I,\fm)$ equals $Z_T(s;\xi',I',\fm)$ if $(\xi, I)$ and
$(\xi', I')$ are equivalent. Let $\fW_\fm$ denote the quotient set of
$W_\ff$ by this equivalence relation. Then for any equivalence
class $\fw\in\fW_\fm$ we can unambiguously define
$Z_T(s;\fw):=Z_T(s;\xi,I,\fm)$ for any $(\xi,I)\in\fw$.
Let $\Clmk$ denote the
ray-class group of $k$ modulo $\fm$. Thus $\Clmk:=\cI_\ff(k)/P_\fm(k)$
where $\cI_\ff(k)$ denotes the
group of fractional ideals prime to $\ff$
and $P_\fm(k)$ the subgroup consisting of those of the form $(a)$ for some
$a\in k^\times,\ a\equiv 1\pmod{\fm}$.
For any $\fc$ in $\Clmk$ and $\fw$ in $\fW_\fm$ we let
$\fc\cdot \fw$ denote the element of $\fW_\fm$ given by
the equivalence class of the pair $(\xi|_{\fa I}, \fa I)\in W_\ff$ where
$(\xi, I)$ is any pair in  the class $\fw$ and
$\fa\in\cI_\ff(k)$ any integral ideal in the class
$\fc$. This map is well-defined and determines
an action of $\Clmk$ on $\fW_\fm$.
One can check that this action is free and transitive.

Let $\cD\lhd\cO$ denote the absolute different of $k$ and write
$\xi_\ff^0$ for the character on $\ff\inv\cD\inv$ which sends $a$ to
$\exp(2\pi i\,\TrkQ(a))$. Thus the pair
$(\xi_\ff^0, \ff\inv\cD\inv)$ lies in $W_\ff$
and we write $\wmz$ for its class in $\fW_\fm$. Let $k(\fm)\subset\bbC$ be
the \emph{ray-class field} over $k$ modulo $\fm$. The Galois
group $G_\fm:=\Gal(k(\fm)/k)$ is isomorphic to $\Clmk$ \emph{via} the
Artin map which sends $\fc\in\Clmk$ to $\sigma_{\fc,\fm}=\sigma_\fc\in G_\fm$.
We let $\bbC G_\fm$ denote the complex group-ring of $G_\fm$ and make the
\begin{defn}\label{defn:1B}
For any cycle $\fm=\ff\fz$ for $k$ and any finite set $T$ of prime
ideals of $\cO$, we write $\PhimT$ for the function
\displaymapdef{\PhimT}{\{s\in\bbC:\Re(s)>1\}}
{\bbC G_\fm}{s}{{\displaystyle\sum_{\fc\in\Clmk}}Z_T(s;\fc\cdot\wmz)\sigma_\fc\inv}
\end{defn}
The basic properties of $\PhimT(s)$ are given in~\cite[\S 3]{twizas}
and \cite[\S 2]{z1}. In particular Theorem~2.2 of \cite{z1} gives a
relation between $\PhimT$ and the primitive $L$-functions of the
characters of $G_\fm$ (or $\Clmk$).

\subsection{$p$-Adic interpolation}\label{subsec:padint}
We turn now to the definition of the $p$-adic analogue of $\Phi_{\fm,T}$
by $p$-adic interpolation. For this we need $k$ to
be \emph{totally real} which we shall assume henceforth.

We choose a prime number $p$, write $\bbC_p$ for the completion of an
algebraic closure of the field $\bbQ_p$ of $p$-adic numbers and fix an
embedding $j:\bar{\bbQ}\rightarrow\bbC_p$.  We let $\mu(\bbQ_p)$ be
the group of roots of unity in $\bbQ_p$,
$w_p$ its cardinality, and
consider the set of rational integers defined by
\[
\cM(p):=
\{
m\in\bbZ\,:\,m\leq 0,\ m\equiv 1 \pmod{w_p}
\}
\]
Let $L(s,\psi)$ denote the
complex $L$-function of a primitive ray-class character
$\psi$. It is well-known that its values at the points of $\cM(p)$ are algebraic
and that their images under $j$ may be `interpolated' to define a $p$-adic $L$-function
attached to the primitive $p$-adic-valued ray-class character $j\circ\psi$
(this is summarised in \cite[Theorem~2.4]{z1}).
It can also be shown (see \cite[Theorem/Definition~2.1]{z1})
that the values of $\PhimT$ on $\cM(p)$ lie in
$\barbbQ G_\fm$ and using, for example,
the $p$-adic $L$-functions one may interpolate these values
whenever the following condition is satisfied
\beql{condit:2A}
\mbox{$T$ contains the set $T_p$ of all primes dividing $p$ in $\cO$.}
\eeq
More precisely, let $D(p)$ denotes
the set $1+2\bbZ_p$, the closure of $\cM(p)$ in $\bbQ_p$.
Then under condition~\ref{condit:2A} there
exists a unique $p$-adic valued function $\PhimTp(s)$, defined on
$D(p)$ and depending on $j$, such that
\beql{eq:1M.2}
\PhimTp(m)=j(\PhimT(m))
\ \ \ \forall
m\in\cM(p).
\eeq
(here, $j$ has been extended to a homomorphism
$\barbbQ G_\fm\rightarrow \bbC_p G_\fm$ by acting on the
coefficients).
Furthermore $\PhimTp$ is $p$-adic meromorphic on $D(p)$. If $\ff$ is
\emph{not} a product of distinct primes lying in $T$, then $\PhimTp$
is actually analytic on this domain. Otherwise it has at most a
unique, simple pole at $s=1$. In all cases $x\PhimTp(s)$ is analytic
in $D(p)$ for any $x$ in the augmentation ideal $I(\bbCp G_\fm)$ of
$\bbCp G_\fm$. Note that we shall write $\Phi_{\fm,T,p}^{(j)}$ instead
of $\PhimTp$ whenever the dependence on $j$ needs to be made explicit.
\subsection{The conjectures}\label{subsec:theconj}
We recall the combined conjectures stated in \cite[\S 3]{z1}.
These are made up of a complex and a $p$-adic part formulated side
by side for the same field $k$ (of degree $r$ say, over $\bbQ$), cycle
$\fm=\ff\fz$ and set $T$, but in terms of $\PhimT(1)$ and $\PhimTp(1)$
respectively. We make the
\begin{hyp}\label{hyp:2A}\ \vspace{-3ex}\\
\begin{enumerate}
\item\label{part:hyp2Apart1}
$k$ is totally real,
\item\label{part:hyp2Apart3}
$\ff$ is not a product of distinct primes lying in $T$
(in particular, $\ff$ is not trivial), and
\item\label{part:hyp2Apart4}
$\fz$ is trivial, \ie\ $\fm=\ff$
\end{enumerate}
\end{hyp}
Hypothesis~\ref{hyp:2A} will be assumed from now on unless the
contrary is explicitly stated, so that, in general, we can write
$\PhifT$ \etc\ in place of $\PhimT$ \etc\ We shall also write $K$ for
the ray-class field $k(\fm)=k(\ff)\subset\bbC$ (necessarily totally
real) and $G$ for $G_\ff=\Gal(K/k)$.  Let $S_\infty$ and
$S_0=S_0(\ff)$ denote respectively the set of infinite (real) places
of $k$ and the set of finite places dividing $\ff$ in $k$. We let
$S=S_\infty\cup S_0$.  The notations $S_\infty(K)$, $S_0(K)$ and
$S(K)$ refer to the sets of places of $K$ dividing those in these
three sets. We abbreviate to $U_S(K)$ or $U_S$ the group $U_{S(K)}(K)$
of $S(K)$-units of $K$.  Let $\iota_1,\ldots,\iota_r$ denote the
embeddings of $k$ into $\barbbQ$ ($\iota_1$ is the inclusion).  For
each $i=1,\ldots,r$ we choose an embedding $\tiota_i:K
\rightarrow\barbbQ$ extending $\iota_i$.  We write $\iota_{i,p}$ for
the $p$-adic embedding $j\circ\iota_i$ of $k$ into $\bbCp$, and also
$\tiota_{i,p}$ for its extension $j\circ\tiota_i:K\rightarrow\bbCp$.
We then define logarithmic maps $\lambda_i:U_S\rightarrow \bbR G$ and
$\lambda_{i,p}:U_S\rightarrow \bbCp G$ by setting
\[
\mbox{
$\lambda_i(u):=
{\displaystyle\sum_{\sigma\in G}}\log|\tiota_i\circ\sigma(u)|\sigma\inv$
and
$\lambda_{i,p}(u):=
{\displaystyle\sum_{\sigma\in G}}
\log_p(\tiota_{i,p}\circ\sigma(u))\sigma\inv$
for all $u\in U_S$.
}
\]
where `$\log_p$' denotes Iwasawa's $p$-adic logarithm.
Both $\lambda_i$ and $\lambda_{i,p}$ are clearly $\bbZ G$-linear
and so `extend' by $\bbQ$-linearity to $\bbQ U_S:=\bbQ\otimes_\bbZ U_S$.
(Henceforth, we shall often write
$\cR A$ in place of $\cR\otimes_\bbZ A$ considered as an $\cR$-module,
for any commutative ring $\cR$ and abelian group $A$.)
These extensions in turn define unique, $\bbQ G$-linear,
group-ring-valued `regulator maps'
$R$ and $R_p$ sending the $r$th exterior power
$\bigwedge_{\bbQ G}^r\bbQ U_S\cong
\bbQ\otimes\bigwedge_{\bbZ G}^r U_S$
into $\bbR G$ and $\bbCp G$ respectively and satisfying
\[
R(u_1\wedge\ldots\wedge u_r)=\det(\lambda_i(u_l))_{i,l=1}^r
\ \ \mbox{and}\ \ \
R_p(u_1\wedge\ldots\wedge u_r)=\det(\lambda_{i,p}(u_l))_{i,l=1}^r
\ \ \forall\,u_1,\ldots,u_r\in \bbQ U_S
\]
We shall use an additive notation for $\WedrQUS$ as $\bbZ G$-module
and write $\lambda_{i,p}^{(j)}$ and $R_p^{(j)}$
instead of $\lambda_{i,p}$
and $R_p$ whenever their dependence on $j$
(\emph{via} the $\tiota_{i,p}$) needs to be made explicit.
Finally, we let $\sqrt{d_k}\in\bbR$ denote the
positive square-root of the (positive) absolute discriminant
$d_k$ of $k$.
\begin{conj}[Basic Combined Conjecture]\label{conj:2C}
If Hypothesis~\ref{hyp:2A} and condition (\ref{condit:2A}) hold then,
there exists $\eta_{\ff,T}\in\bigwedge_{\bbQ G}^r\bbQ U_S$ such that
\beql{eq:2D}
\PhifT(1)=\frac{2^r}{\sqrt{d_k}\prod_{\fp\in T}N\fp}R(\eta_{\ff,T})
\ \ \ \mbox{in $\bbC G$}
\eeq
and
\beql{eq:2E}
\Phi_{\ff,T,p}^{(j)}(1)=
\frac{2^r}{j(\sqrt{d_k})\prod_{\fp\in T}N\fp}R^{(j)}_p(\eta_{\ff,T})
\ \ \ \mbox{in $\bbCp G$}.
\eeq
\end{conj}
\rem\ It is proved in \cite[Prop. 3.3]{z1} that if $\eta_{\ff,T}$
satisfies equation\refeq{2E} for
one embedding $j:\barbbQ\rightarrow\bbC_p$ then it also satisfies it
for any other embedding.\vspace{2ex}\\
In~\cite[\S 3]{z1} `basic' conjectures
were formulated concerning
the existence of elements $\eta$ separately satisfying\refeq{2D}
and\refeq{2E}. These were followed -- under certain
conditions --  by `refined' versions that require $\eta$
to lie in a certain $\bbZ G$-lattice inside
$\bigwedge_{\bbQ G}^r\bbQ U_S$. (These are
`conjectures over $\bbZ$' in the terminology of~\cite{Ta} and~\cite{Ru},
indeed the complex version is closely linked to that of the latter
paper, see~\cite{twizas}, \cite{z1}).
A weakened, combined version of these conjectures was then given
which nevertheless refines Conjecture~\ref{conj:2C}. It is stated
simultaneously for all (eligible) primes $p$.
Before giving this we first recall some notations:
For any $\chi$ in $G^\ast$ (identified with $\Cl_\ff(k)^\ast$) we set
$r(S,\chi):=\rm{dim}_\bbC(e_\chi \bbC U_S)$. Let $\chi_0\in G^\ast$
denote the trivial character and for any place $v$ of $k$, let $G(v)$
denote the decomposition subgroup of $G$ associated to each of the
places $w$ of $K$ dividing $v$.  It can be shown that
\begin{eqnarray}
r(S,\chi)&=&\left\{
\begin{array}{ll}
r+|\{\fq\,:\,\fq|\ff,\ \chi|_{G(\fq)}=1\}|& \mbox{if $\chi\neq\chi_0$, and}\\
r-1+|\{\fq\,:\,\fq|\ff\}|                   & \mbox{if $\chi=\chi_0$}
\end{array}
\right.\nonumber\\
&=&r+\ord_{s=1}(\chi(\PhimT(s)))\label{eq:2Q.2}
\end{eqnarray}
whee the last equation holds
\emph{provided that $(\ff,T)=1$}.
Because $\ff\neq\cO$, it follows that $r(S,\chi)\geq r$ for every
$\chi\in G^\ast$ and $r(S,\chi_0)=r$ if and only if
$\ff$ is a (non-trivial) power of a prime ideal. The latter
condition will be denoted simply `$\ff=\fq^l$'. If it does not hold
then $r(S,\chi_0)>r$ and we shall write `$\ff\neq\fq^l$'. We set
\beql{eq:defids}
e_{S,r}:=\sum_{\chi\in G^\ast\atop r(S,\chi)=r}e_\chi\ \ \
\mbox{and}\ \ \
e_{S,>r}:=1-e_{S,r}=\sum_{\chi\in G^\ast\atop r(S,\chi)>r}e_\chi\ \ \
\eeq
These idempotents actually lie in $\bbQ G$. Let $g$ denote the
cardinality of $G$, then $\tilde{e}_{S,r}:=ge_{S,r}$ and
$\tilde{e}_{S,>r}:=ge_{S,>r}$ clearly lie in $\bbZ G$. For any $\bbZ
G$-module $A$, we shall write $A^{[S,r]}:=\ker \tilde{e}_{S,>r}|A$ so
that $A^{[S,r]}\supset \tilde{e}_{S,r} A\supset gA^{[S,r]}$.
For any $\bbZ G$-submodule $M$ of $U_S$, we denote by
$\overline{\bigwedge^r_{\bbZ G} M}$ the image of the exterior power
$\bigwedge_{\bbZ G}^r M$ in $\WedrQUS$. The conjecture that will be
numerically verified in this article is the following
\begin{conj}[Weak Refined Combined Conjecture]\label{conj:2F}
Suppose that $k$ is totally real and $\ff\neq \cO$
is any proper integral ideal. Then, in the above notations,
there exists a unique element
$\eta_\ff$ of $\left(\WedrQUS\right)\Sr$ with the
following properties
\begin{enumerate}
\item\label{part:conj2E0}
\beql{eq:2S.2}
\frac{2^r}{\sqrt{d_k}}R(\eta_\ff)=\Phi_{\ff,\emptyset}(1)
\eeq
\item\label{part:conj2E1}
For every prime number $p$ with $(p,\ff)=1$ and for every embedding
$j:\barbbQ\rightarrow\bbC_p$ we have
\beql{eq:2R}
\prod_{\fp\in T_p}(1-N\fp\inv\sigma_{\fp,\ff})
\frac{2^r}{j(\sqrt{d_k})}
R_p^{(j)}(\eta_\ff)=\Phi_{\ff,T_p,p}^{(j)}(1)
\eeq
\item\label{part:conj2E2}
If $\ff\neq\fq^l$ then
\beql{eq:2S.5}
\eta_\ff\in\bbZ[1/g]\overline{{\textstyle \bigwedge}^r_{\bbZ G} U_S}\Sr=
\bbZ[1/g]\overline{{\textstyle\bigwedge}^r_{\bbZ G} E(K)}\Sr
\eeq
\item\label{part:conj2E3}
If $\ff=\fq^l$ then
\beql{eq:2T}
\eta_\ff\in\frac{1}{2}\bbZ[1/g]
\overline{{\textstyle \bigwedge}^r_{\bbZ G} U_S}\Sr
\eeq
and
\beql{eq:2U}
I(\bbZ G)\eta_\ff\subset\bbZ[1/g]\overline{{\textstyle \bigwedge}^r_{\bbZ G} U_S}\Sr
\eeq
where $I(\bbZ G)$ is the augmentation ideal of $\bbZ G$.
\end{enumerate}
\end{conj}
\rem\
The point of introducing the
condition $\eta_\ff\in\left(\WedrQUS\right)\Sr$ is that,
essentially without cost,
it allows us to insist upon the uniqueness of $\eta_\ff$
(\cf~\cite[Prop. 3.8]{z1}).
Given this uniqueness and the
relation between $\Phi_{\ff,\emptyset}$ and
$\Phi_{\ff,T_p}$ when $(p,\ff)=1$ (see~\cite[eq.\ (29)]{z1}),
it can be shown that equation\refeq{2R} is actually a consequence of
Conjecture~\ref{conj:2C} (with $T=T_p$) and\refeq{2S.2}.
Moreover, the extra conditions of parts~\ref{part:conj2E2}
and~\ref{part:conj2E3}, which refine Conjecture~\ref{conj:2C},
also follow from it \emph{together with} the assumption of the
`refined complex conjecture' mentioned above for certain
sets $T$. For more details, we refer to Prop.~3.10 of~\cite{z1}.
Note also (\cf\ the preceding remark)
that for given $p$,\refeq{2R} actually holds for all embeddings
$j$ if and only if it holds for one. Finally, for the
(non-conjectural!) equality in\refeq{2S.5}, we refer to~\cite[Lemma 3.5]{z1}.
\section{An expression for $\Phi_{\ff,T_p,p}(1)$ in the quadratic case}\label{quadformula}
\subsection{An application of Shintani's method}
We start with a more general situation than the one suggested by
the title of this section. The data, notated
in the usual way as $k$, $\fm=\ff\fz$, $T$ and $p$, are subject only to
parts~\ref{part:hyp2Apart1} and~\ref{part:hyp2Apart3} of
Hypothesis~\ref{hyp:2A} and to Condition~(\ref{condit:2A}). In
particular, $\ff\neq\cO$.
For any $j:\barbbQ\rightarrow\bbCp$,
we know by Theorem/Definition~2.1 and Lemma~3.3 of~\cite{z1}
that $\PhimTp^{(j)}$ is a $p$-adic analytic map from $D(p)$ to the group ring
$\bbQ_p(\mu_f)^+ G$.
By taking its coefficients we obtain, for each $\fw\in\fW_\fm$,
an analytic map $\ZTpj{\,\cdot\,}{\fw}:D(p)\rightarrow\bbQ_p(\mu_f)^+$
(the \emph{$p$-adic twisted zeta-function attached to $T$, $\fw$ and $j$}).
More precisely, since the action of $\Clmk$ on $\fW_\fm$ is free and
transitive, we can actually {\em define} the $\ZTpj{\,\cdot\,}{\fw}$ by
the equation
\[
\PhimTp^{(j)}(s)=:\sum_{\fc\in\Clmk}\ZTpj{s}{\fc\cdot\wmz}\sigma_\fc\inv
\]
Thus Definition~\ref{defn:1B} and equation\refeq{1M.2}
give the interpolation property
\beql{eq:3A}
\ZTpj{m}{\fw}=j(\ZT{m}{\fw})
\ \ \ \ \mbox{for all $m\in\cM(p)$}
\eeq
which, by density, uniquely characterises $\ZTpj{\,\cdot\,}{\fw}$ as a
continuous map from $D(p)$ to $\bbCp$.

We are interested in calculating $\PhimTp^{(j)}(1)$
in the case $\fz=\emptyset$, $\fm=\ff$ and it clearly suffices to calculate
$\ZTpj{1}{\fc\cdot \fw_\ff^0}$ for each $\fc\in\Clfk$.
However, it is technically and conceptually a little easier to
calculate $\ZTpj{1}{\tilde{\fc}\cdot\fw_{\ff+}^0}$ where the infinite cycle
`$+$' is the formal product of \emph{all}
the real places of $k$ and $\tilde{\fc}$ lies in $\Cl_{\ff+}(k)$. To get back to $\ff$,
we use the natural surjection $\pi_{\ff+,\ff}:\Cl_{\ff+}(k)\rightarrow \Cl_{\ff}(k)$
and the fact that $\ZTpj{1}{\fc\cdot \fw_\ff^0})$ equals
$|\ker\,\pi_{\ff+,\ff}|\,\ZTpj{1}{\tilde{\fc}\cdot\fw_{\ff+}^0}$
for \emph{any} $\tilde{\fc}\in\Cl_{\ff+}(k)$ such that
$\pi_{\ff+,\ff}(\tilde{\fc})=\fc$ (this follows from~\cite[Cor. 2.1]{z1}).
We therefore fix until further notice an element $\fw$ of $\fW_{\ff+}$
and a pair $(\xi,I)\in W_\ff$ representing $\fw$.
Define $f\in\bbZ_{>0}$ by $\ff\cap\bbZ=f\bbZ$ and denote by $\mu_f$ the group
of $f$th roots of unity in $\bbC$. Then
${\rm Im}(\xi)=\mu_f$ (see~\cite[\S 3]{z1}) and when
$\Re(s)>1$, we have
\beql{eq:3B}
\ZT{s}{\fw}=Z_T(s;\xi,I,\ff+)=
\sum_{a\in\cS(I,+,T)/E_{\ff+}}\frac{\xi(a)}{|I:(a)|^s}=
NI^s\sum_{a\in\cS(I,+,T)/E_{\ff+}}\frac{\xi(a)}
{(\iota_1(a)\ldots\iota_r(a))^s}
\eeq
where $\iota_1$,\ldots$,\iota_r\,:\,k\rightarrow\barbbQ\subset\bbC$ are as in
Subsection~\ref{subsec:theconj}.

Shintani's method allows us to analytically continue the second
factor in the fourth member of\refeq{3B} and then find its value
at any $m\in \bbZ_{\leq 0}$. To explain how, we shall simplify matters
by assuming from now on that $k$ is real quadratic ($r=2$). We
require the following notation.
Let $\tau_1$ and $\tau_2$ be two elements of $I\cap\kps$, linearly
independent over $\bbQ$ and such that
$\xi(\tau_1)$, $\xi(\tau_2)\neq 1$.
Then $\tau_1$ and $\tau_2$ define a half-open
`parallelogram' and `cone' in $\kps$ given respectively by:
\[
P(\tau_1,\tau_2):=\{\lambda\tau_1+\mu\tau_2:\lambda,\mu\in\bbQ,\
0<\lambda\leq 1,\ 0\leq\mu<1\}
\]
and
\[
C(\tau_1,\tau_2):=\{\lambda\tau_1+\mu\tau_2:\lambda,\mu\in\bbQ,\
0<\lambda,\ 0\leq\mu\}=
\bigcup_{n_1,n_2\in\bbN}^{\cdot}(P(\tau_1,\tau_2)+n_1\tau_1++n_2\tau_2)
\]
Let $\uiota$
denote the embedding of $k$ into $\bbR^2\cap\barbbQ^2$ which sends
$a\in k$ to $(\iota_1(a),\iota_2(a))$.
Figure~1 illustrates $\uiota (P(\tau_1,\tau_2))$
and $\uiota(C(\tau_1,\tau_2))$ in the case where
$\det\left(
\begin{array}{c}
\uiota(\tau_1)\\
\uiota(\tau_2)
\end{array}
\right)<0$
\begin{figure}
\caption{$\uiota(P(\tau_1,\tau_2))$ and $\uiota(C(\tau_1,\tau_2))$}
\setlength{\unitlength}{4mm}
\[
\begin{picture}(22,17)(-1,-1)
\put(-1,0){\vector(1,0){20}}
\put(0,-1){\vector(0,1){15}}
\thicklines
\put(0,0){\line(1,2){6.5}}
\put(3,6){\line(3,1){9}}
\thinlines
\multiput(0,0)(1.5,0.5){12}{\line(3,1){1}}
\multiput(9,3)(0.75,1.5){4}{\line(1,2){0.5}}
\put(18,-0.3){\makebox(0,0)[t]{$1$}}
\put(-0.3,13){\makebox(0,0)[r]{$2$}}
\put(2.7,6){\makebox(0,0)[r]{$\uiota(\tau_1)$}}
\put(9,2.6){\makebox(0,0)[t]{$\uiota(\tau_2)$}}
\put(0,0){\circle{0.4}}
\put(3,6){\circle*{0.4}}
\put(9,3){\circle{0.4}}
\put(12,9){\circle{0.4}}
\end{picture}
\]
\end{figure}
Clearly, $I\cap P(\tau_1,\tau_2)$ is a fundamental domain for the additive action
of $\bbZ\tau_1+\bbZ\tau_2$ on $I$, so given a class $A\in I/(\bbZ\tau_1+\bbZ\tau_2)$,
we shall write $\tilde{a}=\tilde{a}(A)$ for the unique element of $A\cap P(\tau_1,\tau_2)$.
Then $A=\tilde{a}+\bbZ\tau_1+\bbZ\tau_2$ and
$A\cap C(\tau_1,\tau_2)=\tilde{a}+\bbN\tau_1+\bbN\tau_2$.
We define complex analytic functions
on the set $\{s:\Re(s)>1\}$ (see \eg\ Theorem~\ref{thm:3A} for
convergence and analyticity) by setting
\[
z(s;\xi,A,\tau_1,\tau_2):=\sum_{a\in A\cap C(\tau_1,\tau_2)}
\frac{\xi(a)}{(\iota_1(a)\iota_2(a))^s}=\sum_{n_1,n_2\in\bbN}
\frac{\xi(\tilde{a}+n_1\tau_1+n_2\tau_2)}
{\iota_1(\tilde{a}+n_1\tau_1+n_2\tau_2)^s\iota_2(\tilde{a}+n_1\tau_1+n_2\tau_2)^s}
\]
and also
\[
z(s;\xi,I,\tau_1,\tau_2):=
\sum_{a\in I\cap C(\tau_1,\tau_2)}\frac{\xi(a)}{(\iota_1(a)\iota_2(a))^s}
=
\sum_{A\in I/(\bbZ\tau_1+\bbZ\tau_2)}z(s;\xi,A,\tau_1,\tau_2)
\]
Let $\bbR[[\uX]]$ denote the ring of formal power series
in $\uX=(X_1,X_2)$, a pair of formal variables.
For any pair $\underline{u}=(u_1,u_2)\in\bbR^2$ we write
$(1+\uX)^{\underline{u}}$ for the product $(1+X_1)^{u_1}(1+X_2)^{u_2}$
of two (formal) binomial series in $\bbR[[\uX]]$ and we set
\[
F_{\uiota}(\uX;\xi,A,\tau_1,\tau_2):=
\frac{\xi(\ta)(1+\uX)^{\uiota(\ta)}}
{(1-\xi(\tau_1)(1+\uX)^{\uiota(\tau_1)})(1-\xi(\tau_2)(1+\uX)^{\uiota(\tau_2)})}
\]
and also
\begin{eqnarray}
\lefteqn{F_{\uiota}(\uX;\xi,I,\tau_1,\tau_2):=
\sum_{A\in I/(\bbZ\tau_1+\bbZ\tau_2)}
F_{\uiota}(\uX;\xi,A,\tau_1,\tau_2)=}
\hspace{10em}\nonumber\\
&&\frac{\sum_{\ta\in I\cap P(\tau_1,\tau_2)}\xi(\ta)(1+\uX)^{\uiota(\ta)}}
{(1-\xi(\tau_1)(1+\uX)^{\uiota(\tau_1)})
(1-\xi(\tau_2)(1+\uX)^{\uiota(\tau_2)})}\label{eq:3B.5}
\end{eqnarray}
(The sum in the numerator is of course finite.)
\emph{A priori} these lie
in the fraction field of $\barbbQ[[\uX]]$. However
the constant term $(1-\xi(\tau_1))(1-\xi(\tau_2))$ of their denominators is
non-zero by hypothesis, so they actually lie in $\barbbQ[[\uX]]$ itself (in fact, in
$k(\mu_f)[[\uX]]$, as is easily seen).\vertsp\\
\rem\ Note that intuitively (but illegally) we could also imagine `expanding
the denominator of $F_{\uiota}(\uX;\xi,A,\tau_1,\tau_2)$
as an infinite series in $(1+\uX)$'. We could then write
\[
\mbox{``}\ F_{\uiota}(\uX;\xi,A,\tau_1,\tau_2)=
\sum_{a\in A\cap
C(\tau_1,\tau_2)}\xi(a)(1+\uX)^{\uiota(a)}\ \mbox{''}
\]
We use quotation marks because the sum fails to converge in
$\barbbQ[[\uX]]$, but the idea is useful.\vertsp\\
Let us write $\Delta$ for the differential operator
$(1+X_1)(1+X_2)\frac{\partial^2}{\partial X_1\partial X_2}$
acting on any power series in $X_1$ and $X_2$.
\begin{thm}\label{thm:3A} Let $k$ be a real quadratic field, $\ff\neq\cO$ an integral
ideal with $\ff\cap\bbZ=f\bbZ$, ($f\in\bbZ_{>0}$)
and let $(\xi,I)$ be an element
of $W_\ff$. Suppose that $\tau_1$, $\tau_2$ are two $\bbQ$-linearly
independent
elements of $I\cap k^\times_+\setminus\ker \xi$. Then for any element
$A$ of $I/(\bbZ\tau_1+\bbZ\tau_2)$, we have
\begin{enumerate}
\item\label{part:thm3A1} The function $z(s;\xi,A,\tau_1,\tau_2)$ converges absolutely
for $\Re(s)>1$ and possesses a meromorphic continuation to $\bbC$.
\item\label{part:thm3A2}
For each $m\in\bbZ_{\leq 0}$ this continuation is analytic at $m$ and
for any $\uiota$, we have
\beql{eq:3C}
z(m;\xi,A,\tau_1,\tau_2)=\Delta^{-m}|_{\uX=0}\,F_{\uiota}(\uX;\xi,A,\tau_1,\tau_2)
\eeq
\item\label{part:thm3A3}
$
z(m;\xi,A,\tau_1,\tau_2)\in\bbQ(\mu_f)
$
for all $m\in\bbZ_{\leq 0}$.
\end{enumerate}
\end{thm}
\bPf
Parts~\ref{part:thm3A1} and~\ref{part:thm3A2}
follow from~\cite[Prop. 1]{Shi}
with substitutions ``$r$''=``$n$''=$2$, ``$\chi_i$''$=\xi(\tau_i)$, $i=1,2$, \etc\
Because the ``$\chi_i$'' are different from $1$, the Laurent
series defining Shintani's ``$B_m(a,y,\chi)^{(1)}$'' and
``$B_m(a,y,\chi)^{(2)}$'' are actually power series.
Substituting $1+X_1=e^{-ut_2}$, $1+X_2=e^{-u}$ and
$1+X_1=e^{-u}$, $1+X_2=e^{-ut_1}$ in these two series respectively and
combining them gives\refeq{3C} after a little manipulation.
(Strictly speaking, Shintani's condition that his ``$x_1$''
and ``$x_2$'' be strictly positive is only met if
our $\ta$ lies in the \emph{interior} of $P(\tau_1,\tau_2)$.
However, his proof seems to require only that $\iota(\ta)$
belong to $\bbR^2_+$! In any case, even if $\ta$ \emph{does} lie
on the ray $\bbQ^\times_+\tau_1$, Equation\refeq{3C} can still
be recovered from Shintani's full result (see~\cite[Rem. 2.2]{plim}.))
As for part~\ref{part:thm3A3} of the Theorem, Equation\refeq{3C} already implies that
$z(m;\xi,A,\tau_1,\tau_2)$ lies in $\barbbQ$ (in fact, in $k(\mu_f)$). Now
any $\alpha\in\GQbarQ$ acts coefficientwise on $\barbbQ[[\uX]]$ and
it is clear from the definitions that
\beql{eq:3E}
F_{\uiota}(\uX;\xi,A,\tau_1,\tau_2)^\alpha=
F_{\alpha\circ\uiota}(\uX;\alpha\circ\xi,A,\tau_1,\tau_2)
\eeq
where $\alpha\circ\uiota$ denotes $(\alpha\circ\iota_1,\alpha\circ\iota_2)$.
Since $\Delta$ commutes with $\alpha$, Equations\refeq{3C} and\refeq{3E} give
\begin{eqnarray}
\alpha(z(m;\xi,A,\tau_1,\tau_2))&=&
   \Delta^{-m}|_{\uX=0}\,\left(F_{\uiota}(\uX;\xi,A,\tau_1,\tau_2)^\alpha\right)\nonumber\\
                                &=&
   \Delta^{-m}|_{\uX=0}\,F_{\alpha\circ\uiota}(\uX;\alpha\circ\xi,A,\tau_1,\tau_2)\nonumber\\
                                &=&
   z(m;\alpha\circ\xi,A,\tau_1,\tau_2)\ \ \ \
      \mbox{for all $m\in\bbZ_{\leq 0}$}\label{eq:3E.1}
\end{eqnarray}
and the result follows on letting $\alpha$ run through
$\Gal(\barbbQ/\bbQ(\xi))=\Gal(\barbbQ/\bbQ(\mu_f))$.
\ePf
\rem\ For $\Re(s)>1$, the function $z(s;\xi,A,\tau_1,\tau_2)$ clearly does not
depend on the ordering of the $\iota_i$'s in $\uiota$.
By analytic continuation, neither does the L.H.S.\ of\refeq{3C} and
so the R.H.S.
cannot either. This latter fact was used implicitly in the proof
of part~\ref{part:thm3A3} of the Theorem but is easy to see independently. Indeed if
$\uiota'=(\iota_{\pi(1)},\iota_{\pi(2)})$ for some $\pi\in S_2$
then clearly
$
F_{\uiota'}(X_1,X_2;\xi,A,\tau_1,\tau_2)=
F_{\uiota}(X_{\pi\inv(1)},X_{\pi\inv(2)};\xi,A,\tau_1,\tau_2)
$.
But $\Delta^{-m}|_{\uX=0}$ is symmetric in $X_1$, $X_2$, so
$\Delta^{-m}|_{\uX=0}\,F_{\uiota'}(\uX;\xi,A,\tau_1,\tau_2)=
\Delta^{-m}|_{\uX=0}\,F_{\uiota}(\uX;\xi,A,\tau_1,\tau_2)$, as
required.\vertsp\\
\refstepcounter{example}\noindent{\sc Example \theexample}\label{ex:3A}
By way of illustration, we show how these methods can be used to prove
facts about the values $Z_\emptyset(m;\fw) =
Z_\emptyset(m;\xi,I,\ff+)$ which, in a more general context,
were used in~\cite{z1} (see Lemma 3.2, \emph{ibid.}).
Let $\varepsilon$ be any generator of $E_{\ff+}\cong\bbZ$ and
$\rho$ any element of $I\cap\kps$ not in $\ker\xi$. Then $\varepsilon\rho$
also lies in $I\cap\kps$ but not in $\bbQ\rho$ and
$\xi(\varepsilon\rho)=\xi(\rho)\neq 1$. Thus we can define
$z(s;\xi,I,\rho,\varepsilon\rho)$ and since
it is well known that $I\cap C(\rho,\varepsilon\rho)$
is a fundamental domain for the action of $E_{\ff+}$ on
$I\cap\kps$, Equation\refeq{3B} shows that
\beql{eq:3E.5}
Z_\emptyset(s;\xi,I,\ff+)=
NI^s z(s;\xi,I,\rho,\varepsilon\rho)=
NI^s \sum_{A\in I/(\bbZ\rho+\bbZ\varepsilon\rho)}z(s;\xi,A,\rho,\varepsilon\rho)
\eeq
whenever $\Re(s)>1$. Now apply Theorem~\ref{thm:3A} to equation\refeq{3E.5}.
Part~\ref{part:thm3A1} of the Theorem
allows us to analytically continue the equalities to $s=m\in\bbZ_{\leq 0}$.
Taking $\tau_1=\rho$, $\tau_2=\varepsilon\rho$ in
parts~\ref{part:thm3A2} and~\ref{part:thm3A3}, summing
over $A\in I/(\bbZ\rho+\bbZ\varepsilon\rho)$ and combining
with equation\refeq{3E.5} gives the explicit formula
\[
Z_\emptyset(m;\xi,I,\ff+)=NI^m\Delta^{-m}|_{\uX=0}\,F_{\uiota}(\uX;\xi,I,\rho,\varepsilon\rho)
\ \ \ \ \ \forall\,m\in\bbZ_{\leq 0}
\]
and shows that
\[
Z_\emptyset(m;\xi,I,\ff+)\in\bbQ(\mu_f)
\ \ \ \ \mbox{for all $m\in\bbZ_{\leq 0}$}
\]
The same procedure applied to equation\refeq{3E.1} gives
\[
Z_\emptyset(m;\xi,I,\ff_+)^\alpha=Z_\emptyset(m;\alpha\circ\xi,I,\ff+).
\ \ \ \ \mbox{for all $m\in\bbZ_{\leq 0}$ and $\alpha\in\Gal(\barbbQ/\bbQ)$}
\]
Lemma 3.2 of~\cite{z1} asserts that last two statements hold for
more general $k$ (totally real), $T$ and $\fm$ but the suggested
proof is essentially an elaboration of the above method due to Shintani.
\subsection{Introduction of the prime $p$}
We now introduce a fixed prime number $p$ such that $(p,\ff)=1$.
To $p$-adically interpolate $\Phi_{\fm,T}(s)$ (\ie\
$\ZT{s}{\fw}$) in Subsection~\ref{subsec:padint}, we had to
assume that $T$ contained $T_p$. Therefore, taking $T=T_p$ and
$\tau_1,\tau_2$, to be $\bbQ$-linearly independent and lying in
$I\cap\kps$ but not in $\ker \xi$, as before, we now define a
complex analytic function (see Theorem~\ref{prop:3A}) on the set
$\{s:\Re(s)>1\}$ by setting \beql{eq:3F.5}
z_{T_p}(s;\xi,I,\tau_1,\tau_2):= \sum_{a\in I\cap
C(\tau_1,\tau_2)\atop p\ndiv|I:(a)|}
\frac{\xi(a)}{(\iota_1(a)\iota_2(a))^s} \eeq and also
\[
F_{T_p,\uiota}(\uX;\xi,I,\tau_1,\tau_2):=
\frac{\displaystyle \sum_{\ta\in I\cap P(p\tau_1,p\tau_2)\atop p\ndiv|I:(\ta)|}
\xi(\ta)(1+\uX)^{\uiota(\ta)}}
{(1-\xi(p\tau_1)(1+\uX)^{\uiota(p\tau_1)})
(1-\xi(p\tau_2)(1+\uX)^{\uiota(p\tau_2)})}
\]
which is again an element of $k(\mu_f)[[\uX]]$, since $(p, \ff)=1$ implies
$\xi(p\tau_1),\xi(p\tau_2)\neq 1$.
\begin{thm}\label{prop:3A}
We use the hypotheses and notation of Theorem~\ref{thm:3A}. For any
prime number $p$ with $(p,\ff)=1$ we have:
\begin{enumerate}
\item
The function $z_{T_p}(s;\xi,I,\tau_1,\tau_2)$ converges absolutely
for $\Re(s)>1$ and possesses a meromorphic continuation to $\bbC$.
\item\label{part:prop3A2}
For each $m\in\bbZ_{\leq 0}$ this continuation is analytic at $m$ and
for any $\uiota$, we have
\[
z_{T_p}(m;\xi,I,\tau_1,\tau_2)=\Delta^{-m}|_{\uX=0}\,
F_{T_p,\uiota}(\uX;\xi,I,\tau_1,\tau_2)
\]
\item\label{part:prop3A3}
$z_{T_p}(m;\xi,I,\tau_1,\tau_2)\in\bbQ(\mu_f)$
for all $m\in\bbZ_{\leq 0}$.
\end{enumerate}
\end{thm}
\bPf\ The condition $p\ndiv|I:(a)|$ is equivalent to $a\nin\fp I$ for any prime ideal
$\fp$ of $\cO$ dividing $p$.
Since $\bbZ p\tau_1+\bbZ p\tau_2\subset pI\subset\bigcap_{\fp |p}\fp I$,
it follows that the set of $a$ satisfying this condition is a
union of those cosets $A\in I/(\bbZ p\tau_1+\bbZ p\tau_2)$
not contained in (\ie\ not intersecting) $\fp I$ for any $\fp|p$.
Letting $\cA'$ denote the (finite) set of all such cosets,
it follows from the definitions that
\[
z_{T_p}(s;\xi,I,p\tau_1,p\tau_2)=\sum_{A\in\cA'}z(s;\xi,A,p\tau_1,p\tau_2)
\]
and
\[
F_{T_p,\uiota}(\uX;\xi,I,\tau_1,\tau_2)=\sum_{A\in\cA'}
F_{\uiota}(\uX;\xi,A,p\tau_1,p\tau_2)
\]
But $C(\tau_1,\tau_2)=C(p\tau_1,p\tau_2)$, so
$z_{T_p}(s;\xi,I,p\tau_1,p\tau_2)=z_{T_p}(s;\xi,I,\tau_1,\tau_2)$. The
Proposition therefore follows from Theorem~\ref{thm:3A} with
$p\tau_1$ and $p\tau_2$ in place of $\tau_1$ and $\tau_2$.
\ePf
Before proceeding with a $p$-adic
interpolation of $z_{T_p}(s;\xi,I,\tau_1,\tau_2)$,
we formulate a hypothesis and a definition and prove a lemma.
\begin{hyp}\label{hyp:3A}\ \vspace{-3ex}\\
\begin{enumerate}
\item\label{part:hyp3A1} $(p,\ff)=1$,
\item\label{part:hyp3A2} $p$ splits in $k$, \ie\ $p\cO=\fp_1\fp_2$
with $\fp_1\neq\fp_2$, and
\item\label{part:hyp3A3} $I$ is prime to $p$, \ie\ $\ord_{\fp_i}(I)=0$
for $i=1,2$.
\end{enumerate}
\end{hyp}
\rem\ Condition~\ref{part:hyp3A1} has already been imposed.
Without it the nature of the interpolation problem would change
significantly. Assuming it, and taking $T=T_p$,
Condition~\ref{part:hyp2Apart3} of Hypothesis~\ref{hyp:2A}
is \emph{equivalent to} the condition
$\ff\neq\cO$. Condition~\ref{part:hyp3A2} of Hypothesis~\ref{hyp:3A}
is not necessary for (a generalized version of) the results that follow
but it simplifies their
exposition and the computations based on them.
Finally, Condition~\ref{part:hyp3A3} is no obstruction to
computing $Z_{T_p}(s;\fw)$ since any $\fw\in\fW_{\ff+}$
can always be represented by some $(\xi,I)$ with $I$ prime to $p$.
\begin{defn} Let $\{c_n\}_{n=1}^\infty$ be the sequence of
rational integers given by
\beql{eq:3H}
c_n=c_n(p):=\sum_{\zeta^p=1}(\zeta-1)^n=
p\sum_{0\leq r\leq n/p}(-1)^{n-pr}\left(
\begin{array}{c}
n\\
pr
\end{array}
\right)
\eeq
where $\zeta$ runs through the $p$th roots of unity
in any algebraic closure of $\bbQ$.
\end{defn}
(The second formula in\refeq{3H} follows from the first by expanding
$(\zeta-1)^n$.)
Let $P(X)$ denote the polynomial $((X+1)^p-1)/X$.
Taking $n\geq p$, writing out the L.H.S. of the equation
$(\zeta-1)^{n-(p-1)}P(\zeta-1)=0$ as a polynomial in $\zeta-1$
and summing over $\zeta$, we obtain the useful recurrence relation
\beql{eq:3I}
c_n=-\left(
\left(
\begin{array}{c}
p\\
1
\end{array}
\right)c_{n-1}+
\left(
\begin{array}{c}
p\\
2
\end{array}
\right)c_{n-2}+
\ldots
+\left(
\begin{array}{c}
p\\
p-1
\end{array}
\right)c_{n-(p-1)}
\right)\ \ \ \forall\,n\geq p
\eeq
with the initial conditions $c_n=(-1)^np$ for $1\leq n\leq (p-1)$
which follow from the second formula in\refeq{3H}.
Let $|\cdot|_p$ denote the absolute value
on $\bbC_p$ normalised by $|p|_p=p\inv$
and for $x\in\bbR$, let $\lceil x \rceil$
denote ${\rm min}\{l\in\bbZ\,:\,x\leq l\}$.
\begin{lemma}\label{lemma:estimcn} For all $n\geq 1$, we have
$|c_n|_p\leq p^{-\lceil n/(p-1)\rceil}$ and the quotient
$c_n/pn$ is $p$-integral.
\end{lemma}
\bPf
The estimate follows from the fact that $|\zeta-1|_p=p^{-1/(p-1)}$ for
every $\zeta$ not equal to $1$, or by induction from\refeq{3I}. For the
$p$-integrality statement, define $m\in\bbN$ by $p^m\leq n<p^{m+1}$
and note that $\ord_p(c_n/pn)$ is at least
$\lceil n/(p-1)\rceil-1-\ord_p(n)$ which is clearly zero if
$m=0$ and is otherwise at least
$\lceil p^m/(p-1)\rceil-1-m=\sum_{i=0}^{m-1}(p^i-1)\geq 0$.\ePf
We can now state the main result of this section.
\begin{thm}\label{thm:3B}
We use the hypotheses and notation of Theorem~\ref{thm:3A}. Suppose that
$p$ is any prime number satisfying Hypothesis~\ref{hyp:3A} and
$j:\barbbQ\rightarrow\bbCp$ any embedding. We have:
\begin{enumerate}
\item\label{part:thm3B1}
There exists a unique $p$-adically continuous function
$z_{T_p,p}^{(j)}(\,\cdot\,;\xi,I,\tau_1,\tau_2)\,:\,D(p)\longrightarrow \bbCp$
satisfying the interpolation condition
\beql{eq:3I.1}
z_{T_p,p}^{(j)}(m;\xi,I,\tau_1,\tau_2)=j(z_{T_p}(m;\xi,I,\tau_1,\tau_2))
\ \ \ \forall
m\in\cM(p)
\eeq
\item\label{part:thm3B2}
For any $\uiota$, write
$(1+\uX)\inv F_{\uiota}(\uX;\xi,I,\tau_1,\tau_2)^j$ as $\sum_{i,l\geq 0}a_{i,l}X_1^iX_2^l$.
Then $a_{i,l}$ lies in $\bbZp[\mu_f]$ for all $i,l\geq 0$ and
\beql{eq:3I.2}
z_{T_p,p}^{(j)}(1;\xi,I,\tau_1,\tau_2)=
\frac{1}{p^2}\sum_{i,l\geq 0}\frac{c_{i+1}c_{l+1}a_{i,l}}{(i+1)(l+1)}\in\bbZp[\mu_f]
\eeq
\end{enumerate}
\ePf
\end{thm}
\emph{Note that we do mean $F_{\uiota}(\uX;\xi,I,\tau_1,\tau_2)^j$,
not} $F_{T_p,\uiota}(\uX;\xi,I,\tau_1,\tau_2)^j$ in part~\ref{part:thm3B2} and that
the exponent indicates that $j$ has been applied to the coefficients of the power-series.
The estimate of $\ord_p(c_n/pn)$ in the proof of Lemma~\ref{lemma:estimcn}
shows that the infinite sum in\refeq{3I.2} converges ($p$-adically).
The proof of Theorem~\ref{thm:3B}
will provide an expression for $z_{T_p,p}^{(j)}(s;\xi,I,\tau_1,\tau_2)$
as a $p$-adic integral (see Equation\refeq{3P.5}).
We defer it while we deduce a result that allows us in principle to
calculate $Z^{(j)}_{T_p,p}(1;\fw)$ for all $\fw\in\fW_{\ff+}$ and hence
$\Phi_{\ff,T_p,p}(1)$ as explained at the beginning of this section.
Let $\omega:\bbZ_p^\times\rightarrow\mu(\bbQ_p)$ be the Teichm\"{u}ller character
which is uniquely defined by the requirement that $\langle x\rangle:=\omega\inv(x)x$
should lie in $1+p\bbZ_p$ for all $x\in\bbZ_p^\times$ (and in $1+4\bbZ_2$ if $p=2$).
\begin{cor}\label{cor:3B} Under the hypotheses of the Theorem, let
$\fw$ be the class of $(\xi,I)$ in $\fW_{\ff+}$,
let $\varepsilon$ be any generator of
$E_{\ff+}$ and let $\rho$ any element of $I\cap\kps$ not in $\ker\xi$.
Then
\beql{eq:3I.5}
Z^{(j)}_{T_p,p}(s;\fw)=\omega(NI)\langle NI\rangle^s
z_{T_p,p}^{(j)}(s;\xi,I,\rho,\varepsilon\rho)\ \ \ \ \forall\,s\in D(p)
\eeq
where $z_{T_p,p}^{(j)}(s;\xi,I,\rho,\varepsilon\rho)$ is as in
part~\ref{part:thm3B1} of the Theorem. In particular,
\beql{eq:3J}
Z^{(j)}_{T_p,p}(1;\fw)=
\frac{NI}{p^2}\sum_{i,l\geq 0}\frac{c_{i+1}c_{l+1}a_{i,l}}{(i+1)(l+1)}\in\bbZp[\mu_f]
\eeq
where $\fw$ is the class of $(\xi,I)$ in $\fW_{\ff+}$ and the $a_{i,l}$
are defined by
$(1+\uX)\inv
F_{\uiota}(\uX;\xi,I,\rho,\varepsilon\rho)^j=\sum_{i,l\geq 0}a_{i,l}X_1^iX_2^l$
for any $\uiota$.
\end{cor}
Note that Hypothesis~\ref{hyp:3A}~\ref{part:hyp3A3} implies
$(NI,p)=1$, so that $\omega(NI)$ is well-defined and the function
$\langle NI\rangle^s$ is both well-defined and analytic for $s\in\bbZ_p$.\vertsp\\
\bPf\ Arguing just as in Example~\ref{ex:3A}, Equations\refeq{3B}
and\refeq{3F.5} show that $Z_{T_p}(s;\fw)$ equals
$NI^s z_{T_p}(s;\xi,I,\rho,\varepsilon\rho)$ for every $s$ such that
$\Re(s)>1$ and hence, by analytic continuation, for every $m\in\cM(p)$. But
$NI^m=\omega (NI)\langle NI \rangle^m$ for such $m$,
so part~\ref{part:thm3B1} of the Theorem gives
\[
\omega (NI)\langle NI \rangle^m
z_{T_p,p}^{(j)}(m;\xi,I,\rho,\varepsilon\rho)=
NI^m j(z_{T_p}(m;\xi,I,\rho,\varepsilon\rho))=j(Z_{T_p}(m,\fw))
\ \ \forall\,m\in\cM(p)
\]
Since the function $\omega(NI)\langle NI\rangle^s
z_{T_p,p}^{(j)}(s;\xi,I,\rho,\varepsilon\rho)$ is continuous on
$D(p)$ by part~\ref{part:thm3B1} of the Theorem,
the equality\refeq{3I.5} now follows from the uniqueness of the interpolation
in\refeq{3A}. Equation\refeq{3J} is then a direct consequence of
part~\ref{part:thm3B2} of the Theorem.\ePf\
\subsection{Proof of theorem~\ref{thm:3B}}\label{subsec:pfthm3B}
Our methods of proof generalise those
employed by Lang to evaluate $L_p(1,\chi)$ in~\cite[Ch. 4]{La}.
They use the theory of $p$-adic-valued measures on $\bbZp^2$
and their relation to formal power series. We recall that
such a measure is a $\bbC_p$-valued, bounded linear functional on the
$\bbC_p$-Banach algebra $\ContZps$ of all continuous,
$\bbC_p$-valued functions on $\bbZ_p^2$ under the (ultrametric)
uniform norm $\|\cdot\|$.
The `boundedness' requirement on such a functional $\nu$
means that the set
$\{|\nu(f)|_p/\|f\|\,:\,f\in\ContZps, f\neq 0\}$ is bounded.
By writing $\|\nu\|$ for its supremum we define an ultrametric
norm under which the set $\MeasZps$ of all such
measures acquires the structure of a $\bbC_p$-Banach space.
For $\nu\in\MeasZps$ and $f\in\ContZps$ the value $\nu(f)$,
will often be written as $\int_{\bt\in\bbZp^2}f(\bt)\,d\nu$ or just $\int f\,d\nu$.
Clearly, $\MeasZps$ is a natural
$\ContZps$-module where for $\nu$ in the former
and $g$ in the latter, we define the measure $g\nu$ by
$
\int f\,d(g\nu):=\int fg\,d\nu\ \ \forall\,f\in\ContZps.
$
When $g$ is the characteristic function $\chi_S$ of
an open and closed subset $S$ of $\bbZp^2$, we often write
$\nu|_S$ for $\chi_S\nu$ (`the restriction of $\nu$ to $S$') and
$\int_{\bt\in S}f(\bt)\,d\nu$ instead of
$\int f\,d(\chi_S\nu)=\int f\chi_S\,d\nu$.

Let us write $\cA(\uX)$ for the
$\bbCp$-subspace of $\bbC_p[[\uX]]$ consisting of those
power-series with ($p$-adically) bounded coefficients.
For such a power series $F$, we define the norm $\|F\|$
to be the supremum of the $p$-adic absolute values of its
coefficients. The power-series/measure correspondence is then a
norm-preserving isomorphism of $\bbCp$-Banach spaces between
$\cA(\uX)$ and $\MeasZps$. In~\cite[Ch. 4]{La},
Lang discusses in detail a restricted correspondence
${\mathbb O}[[X]]\leftrightarrow {\rm Meas}^{(1)}(\bbZp,\bbC_p)$
where $\mathbb O$ denotes the ring $\{a\in\bbCp:|a|_p\leq 1\}$
and ${\rm Meas}^{(1)}(\bbZp,\bbC_p)$ the space of
measures on $\bbZp$ of norm $\leq 1$. By simply taking $\bbCp$-spans we get a
correspondence between $\cA(X)$ and ${\rm Meas}(\bbZp,\bbC_p)$
(see also Appendices~5 and~6 of~\cite{Sch}).
The generalisation of this from one to two (or more) variables
seems to be well known although we have been unable to find a
full and detailed account in the published literature.
In any case, it is very straightforward. The facts we require are
as follows (see also~\cite{alshig2} for the $r$-variable case).
The correspondence can be characterised as associating
$\nu\in\MeasZps$ with $F\in\cA(\uX)$ if and only if
\beql{eq:3L}
\int\limits_{\bt\in\bbZp^2}(1+u_1)^{t_1}(1+u_2)^{t_2}
\,d\nu=F(u_1,u_2)
\ \ \ \ \ \ \mbox{$\forall\,(u_1,u_2)\in\bbCp^2$ s.t.
$|u_1|_p$, $|u_2|_p<1$}
\eeq
in which case we shall write $\nu=\cN(F)$ and
$F=\cF(\nu)$.
By expanding Equation\refeq{3L} as a power series in
$u_1$ and $u_2$ we can deduce (see~\cite{alshig2}) that for all
$\nu\in\MeasZps$ and $n_1,n_2\in\bbZ_{\geq 0}$
\beql{eq:3M}
\cF(t_1^{n_1}t_2^{n_2}\nu)=
\left((1+X_1)\frac{\partial}{\partial X_1}\right)^{n_1}
\left((1+X_1)\frac{\partial}{\partial X_1}\right)^{n_2}\cF(\nu)
\eeq
We denote by $D^-$ the open $p$-adic bidisc
$\{(x_1,x_2)\in\bbCp^2\,:\,|x_1|_p,\ |x_2|_p<1\}$
and by $\cA_1(\uX)\subset\bbCp[\uX]$ the $\bbCp$-algebra of power series
convergent at every point of $D^-$. Thus
$\cA_1(\uX)$ contains $\cA(\uX)$.
We define an action of the group $\mu_p^2=\mu_p(\bbCp)^2$
on $\cA_1(\uX)$ by setting
$
(\underline{\zeta}\bullet F)(X_1,X_2)=
F(\zeta_1(1+X_1)-1,\zeta_2(1+X_2)-1)
$
for any $\underline{\zeta}=(\zeta_1,\zeta_2)\in \mu_p^2$
and $F\in\cA_1(\uX)$. It is easy to check that this indeed gives
a well-defined $\bbCp$-linear left action (see~\cite[Sec. 3.2]{plim}).
Restricting to the subgroup $\mu_p\times\{1\}\subset\mu_p^2$, the
idempotent corresponding to the trivial character of this
group is the operator $V_1$:
\[
V_1\bullet F:=\frac{1}{p}\sum_{\zeta_1^p=1}F(\zeta_1(1+X_1)-1,X_2)
\]
An operator $V_2$ is defined similarly by acting on the
variable $X_2$ and we write $U$ for the idempotent operator
$U=(1-V_1)(1-V_2)=(1-V_2)(1-V_1)$.
It can be checked (\cite[Sec. 3.2]{plim})
that the `$\bullet$' action preserves $\cA(\uX)$. Moreover,
it follows easily from\refeq{3L} that given any $F\in\cA(\uX)$ and
$(\zeta_1,\zeta_2)\in \mu_p^2$, the measure
$\cN((\zeta_1,\zeta_2)\bullet F)$ is simply the measure
$\cN(F)$ multiplied by the continuous (locally constant) function
$\bt\mapsto\zeta_1^{t_1}\zeta_2^{t_2}$. From this it follows that
$
\cN(V_1\bullet F)=\chi_{p\bbZp\times\bbZp}\cN(F)
$,
$
\cN(V_2\bullet F)=\chi_{\bbZp\times p\bbZp}\cN(F)
$
and so
\beql{eq:3M.5}
\cN(U\bullet F)=\chi_{(\bbZp^\times)^2}\cN(F)
\eeq
We fix $\uiota$, $\tau_1$, $\tau_2$ and $j$,
and abbreviate the power series $F_{\uiota}(\uX,\xi,\tau_1,\tau_2)^j$ and
$F_{T_p, \uiota}(\uX,\xi,\tau_1,\tau_2)^j$ to
$F_\xi$ and $F_\xi^\ast$ respectively.
\begin{lemma}\label{lemma:3A}
$F_\xi$ lies in $\bbZp[\mu_f][[\uX]]$, hence in
$\cA(\uX)$. Furthermore $F_\xi^\ast=U\bullet F_\xi$.
\end{lemma}
\bPf\
On the R.H.S of\refeq{3B.5} we can multiply both the
numerator and denominator by the power series
\begin{eqnarray*}
\lefteqn{
\left(\sum_{n_1=0}^{p-1}(\xi(\tau_1)(1+\uX)^{\uiota(\tau_1)})^{n_1}\right)
\left(\sum_{n_2=0}^{p-1}(\xi(\tau_2)(1+\uX)^{\uiota(\tau_2)})^{n_2}\right)=}
\hspace{18em}\\
&&\sum_{n_1,n_2=0}^{p-1}\xi(n_1\tau_1+n_2\tau_2)
(1+\uX)^{\uiota(n_1\tau_1+n_2\tau_2)}
\end{eqnarray*}
But $I\cap P(p\tau_1,p\tau_2)$ is the disjoint union of the
translates $n_1\tau_1+n_2\tau_2+(I\cap P(\tau_1,\tau_2))$ for $0\leq n_1,n_2\leq p-1$,
so (after applying $j$ to\refeq{3B.5}) we see that
$F_\xi$ can be written as
\beql{eq:3N}
F_{\xi}=
\sum_{\ta\in I\cap P(p\tau_1,p\tau_2)}
\frac{\xi(\ta)(1+\uX)^{\uiota_p(\ta)}}
{(1-\xi(p\tau_1)(1+\uX)^{\uiota_p(p\tau_1)})
(1-\xi(p\tau_2)(1+\uX)^{\uiota_p(p\tau_2)})}
\eeq
Here for any $a\in k$, the notation $(1+\uX)^{\uiota_p(a)}$
indicates $((1+\uX)^{\uiota(a)})^j=(1+X_1)^{a_1}(1+X_2)^{a_2}$, the product
of two formal \emph{$p$-adic} binomial series with
$a_1:=j\circ\iota_1(a)$ and $a_2:=j\circ\iota_2(a)$ which are the two
embeddings of $a$ in $\bbCp$. Now, parts~\ref{part:hyp3A2}
and~\ref{part:hyp3A3} of Hypothesis~\ref{hyp:3A} imply that
$a_1$ and $a_2$ lie in $\bbZp$ whenever $a$ lies in $I$.
As is well known, this implies in turn that the series $(1+X_1)^{a_1}$
and $(1+X_2)^{a_2}$ have coefficients in $\bbZp$ for any such $a$,
hence that the numerator and the (common)
denominator of each term on the R.H.S. of\refeq{3N}
lie in $\bbZp[\mu_f][[\uX]]$. The constant term of this
denominator is $(1-\xi(\tau_1)^p)(1-\xi(\tau_2)^p)=:c$, say. Now,
Hypothesis~\ref{hyp:3A}~\ref{part:hyp3A1} of  implies that
$(p,f)=1$ so $\xi(\tau_1)$ and $\xi(\tau_2)$ are roots of unity of
order prime to $p$, non-trivial by assumption, so the same is true
of their $p$th powers. It follows that $c$
lies in $\bbZp[\mu_f]^\times$, so that each term in\refeq{3N}
actually has denominator lying
in $\bbZp[\mu_f][[\uX]]^\times$ and hence itself lies in
$\bbZp[\mu_f][[\uX]]\subset\cA[[\uX]]$. The first statement in the
Lemma follows. As for the second, it is easy to show that
for any $a\in I$, the element $(\zeta_1,\zeta_2)\in\mu_p^2$
acts on $(1+X_1)^{a_1}(1+X_2)^{a_2}$ by multiplication by
$\zeta_1^{a_1}\zeta_2^{a_2}$ (since $a_1,a_2\in\bbZp$).
Thus if $\ta$ is an element of $I\cap P(p\tau_1,p\tau_2)$
then $(\zeta_1,\zeta_2)\in\mu_p^2$
multiplies the corresponding term on the R.H.S. of\refeq{3N}
by $\zeta_1^{\ta_1}\zeta_2^{\ta_2}$
(using the fact that it acts trivially on the denominator).
It follows that $V_1$ acts on this term by $1$ or $0$ respectively,
according as $p$ does or does not divide $\ta_1$ in $\bbZp$, and similarly for $V_2$,
\emph{mutatis mutandi} with the result that $U$ acts by $0$ or $1$
according as $p$ does or does not divide $\ta_1\ta_2=|\cO:I||I:(\ta)|$.
Here $|\cO:I|\in\bbQ^\times$ is the (generalised) index and lies in $\bbZp^\times$ by
Hypothesis~\ref{hyp:3A}~\ref{part:hyp3A3}. Putting this all together and applying
$U$ to Equation\refeq{3N} we get
\[
U\bullet F_{\xi}=
\sum_{\ta\in I\cap P(p\tau_1,p\tau_2)\atop p\ndiv|I:(\ta)|}
\frac{\xi(\ta)(1+\uX)^{\uiota_p(\ta)}}
{(1-\xi(p\tau_1)(1+\uX)^{\uiota_p(p\tau_1)})
(1-\xi(p\tau_2)(1+\uX)^{\uiota_p(p\tau_2)})}
\]
and the R.H.S. is, by definition, the image of
$F_{T_p,\uiota}(\uX;\xi,I,\tau_1,\tau_2)$ under $j$, as required.
\ePf
Lemma~\ref{lemma:3A} implies that both $F_\xi$ and $F_\xi^\ast$
lie in $\cA(\uX)$ and that if we set
$\nu_\xi=\cN(F_\xi)$ and $\nu_\xi^\ast=\cN(F_\xi^\ast)$
then $\nu_\xi^\ast=\chi_{(\bbZp^\times)^2}\nu_\xi$.
For any elements $m$ of $\cM(p)$ and $F$ of
$\cA(\uX)$, Equation\refeq{3M} implies that $\cN(\Delta^{-m}F)=(t_1t_2)^{-m}\cN(F)$
so Equation\refeq{3L} with $u_1=u_2=0$ gives
\[
\int\limits_{\bbZp^2}(t_1t_2)^{-m}d\cN(F)=\Delta^{-m}|_{\uX=0}F
\]
Applying this with $F=F_\xi^\ast$, and noting that $\Delta$
commutes with $j$, Theorem~\ref{prop:3A} part~\ref{part:prop3A2} gives,
for all $m\in\cM(p)$:
\begin{eqnarray}
j(z_{T_p}(m;\xi,I,\tau_1,\tau_2))&=&\Delta^{-m}|_{\uX=0}F_\xi^\ast
                   =\int\limits_{\bbZp^2}(t_1t_2)^{-m}d\nu_\xi^\ast\nonumber\\
                                &=&\int\limits_{(\bbZp^\times)^2}(t_1t_2)^{-m}d\nu_\xi
                   =\int\limits_{(\bbZp^\times)^2}
                        \omega(t_1t_2)\inv\langle t_1t_2\rangle^{-m}d\nu_\xi\label{eq:3P}
\end{eqnarray}
For any $s\in D(p)$ we define
\displaymapdef{f_s}{\bbZp^2}{\bbZp}{\bt}
{
\left\{
\begin{array}{ll}
\omega(t_1t_2)\inv\langle t_1t_2\rangle^{-s}& \mbox{if $\bt\in (\bbZp^\times)^2$}\\
0& \mbox{otherwise}
\end{array}
\right.
}
Thus with our definitions, the last integral in\refeq{3P} is strictly to be
interpreted as $\int_{\bbZp^2}f_m(\bt)d\nu_\xi$. But
$f_s(\bt)$ is easily seen to be uniformly continuous as a function of
$(s,\bt)\in D(p)\times\bbZp^2$, so it follows that on defining
\beql{eq:3P.5}
z_{T_p,p}^{(j)}(s;\xi,I,\tau_1,\tau_2):=
\int\limits_{\bbZp^2}f_s(\bt)d\nu_\xi=
\int\limits_{(\bbZp^\times)^2}\omega(t_1t_2)\inv\langle t_1t_2\rangle^{-s}d\nu_\xi
\eeq
we have a $p$-adically continuous function which, by\refeq{3P}, satisfies the interpolation
condition\refeq{3I.1}. Its unicity follows from the density of $\cM(p)$
in $D(p)$. This proves part~\ref{part:thm3B1} of
Theorem~\ref{thm:3B} and gives
\beql{eq:3Q}
z_{T_p,p}^{(j)}(1;\xi,I,\tau_1,\tau_2)=\int\limits_{(\bbZp^\times)^2}(t_1t_2)^{-1}d\nu_\xi=
\int\limits_{\bbZp^2}f_1(\bt)d\nu_\xi=G_\xi(0,0)
\eeq
where we define $G_\xi(\uX)\in\cA(\uX)$ to be $\cF(f_1\nu_\xi)\in\cA(\uX)$.
To determine $G_\xi$ we use the
\begin{lemma}\label{lemma:3B}
If $H$ is any element of $\cA_1(\uX)$ satisfying
$\Delta H=F_\xi$ then $U\bullet H=G_\xi$.
(In particular, $U\bullet H$ lies in $\cA(\uX)$.)
\end{lemma}
\bPf\ Since $\Delta$ commutes with the $\bullet$-action,
the condition on $H$, together with Equations\refeq{3M} and\refeq{3M.5},
implies that
\[
\Delta(U\bullet H-G_\xi)=\Delta U\bullet H-\cF(t_1t_2f_1(\bt)\nu_\xi)=
U\bullet(\Delta H)-\cF(\chi_{(\bbZp^\times)^2}\nu_\xi)=
U\bullet F_\xi-U\bullet F_\xi=0
\]
Since $(1+\uX)$ is an invertible power series,
it follows from the definition of $\Delta$ that $U\bullet H-G_\xi=B_1(X_1)+B_2(X_2)$
for some single-variable power series $B_1$
and $B_2$. Since $U\bullet H-G_\xi$ lies in $\cA_1(\uX)$, it is easy
to see that both $B_1(X_1)$ and $B_2(X_2)$ must too, and also
that $V_2\bullet B_1(X_1)=B_1(X_1)$ and $V_1\bullet
B_2(X_2)=B_2(X_2)$. Thus $U\bullet B_1(X_1)=U\bullet B_2(X_2)=0$.
On the other hand, $U\bullet G_\xi=G_\xi$ (since
$\chi_{(\bbZp^\times)^2}f_1=f_1$) and since $U$ is idempotent, we
obtain
\[
U\bullet H-G_\xi=U\bullet(U\bullet
H-G_\xi)=U\bullet(B_1(X_1)+B_2(X_2))=0
\]
proving the Lemma.\ePf
Now let us write $(1+\uX)\inv F_\xi=\sum_{i,l\geq 0}a_{i,l}X_1^iX_2^l$
as in the statement of the Theorem. Lemma~\ref{lemma:3A}
implies that the $a_{i,l}$ lie in $\bbZ_p[\mu_f]$ and so, by
easy, standard estimates, the power series $H_0$ defined by
\[
H_0(\uX):=\sum_{i,l\geq 0}\frac{a_{i,l}}{(i+1)(l+1)}X_1^{i+1}X_2^{l+1}
\]
lies in $\cA_1(\uX)$. Moreover, $\Delta H_0$ equals
$F_\xi$ by construction. Therefore, Lemma~\ref{lemma:3B} gives
\beql{eq:3Q.5}
G_\xi=U\bullet H_0
\eeq
Now, for any $\zeta\in\mu_p$ we clearly have $H_0(\zeta-1,0)=H_0(0,\zeta-1)=H_0(0,0)=0$
from which it follows in particular that
$(V_1\bullet H_0)(0,0)=(V_2\bullet H_0)(0,0)=0$. Thus,
combining Equations\refeq{3Q} and\refeq{3Q.5},
and expanding $U$ as $1-V_1-V_2+V_1V_2$ we obtain
\begin{eqnarray}
\lefteqn{z_{T_p,p}^{(j)}(1;\xi,I,\tau_1,\tau_2)=(U\bullet H_0)(0,0)=
(V_1V_2\bullet H_0)(0,0)=}\hspace{10em}\nonumber\\
&&
\frac{1}{p^2}\sum_{\zeta_1^p=\zeta_2^p=1}H_0(\zeta_1-1,\zeta_2-1)=
\frac{1}{p^2}\sum_{i,l\geq 0}\frac{c_{i+1}c_{l+1}a_{i,l}}{(i+1)(l+1)}
\label{eq:3R}
\end{eqnarray}
Finally, Lemma~\ref{lemma:estimcn} shows that the last member of\refeq{3R}
lies in $\bbZp[\mu_f]$.\ePf
\section{Numerical investigation of conjecture \ref{conj:2F}}\label{vericonj}
In this section we present a number of examples in which
Conjecture~\ref{conj:2F} is verified up to the precision of computation
for a real quadratic field $k$.
Of course, in each example,
part~\ref{part:conj2E1} of the conjecture will only be checked for a (small)
finite set of primes! (Moreover, these primes will be subjected
to certain further conditions that facilitate
the calculation of $\PhifTpp(1)$. See below.)
Before presenting the examples themselves, we explain some of
our computational techniques and methods.
\subsection{Remarks on computational methods}\label{comp}
To compute $\Phi_{\ff,\emptyset}(1)$, we use the decomposition given
by \cite[eq.\ (12)]{z1} and take $s = 1$ in \cite[eq.\ (10)]{z1},
making use of the functional equation of $L(s,\tchi)$, to obtain the
expression
$$
\Phi_{\ff,\emptyset}(1) = \frac{4}{\sqrt{d_k}}
\sum_{\chi \in G^\ast \atop \chi \not= \chi_0}
\prod_{\fp\ \mathit{prime} \atop \fp|\ff, \fp\ndiv\ff(\chi)}
(1 - \tchi(\fp)\inv) L^{(2)}(0, \tchi\inv) e_\chi +
\left\{
\begin{array}{ll}
0 & \mbox{if $\ff \not= \fq^l$,} \\
- \log(N\fq) \frac{2 h_k R_k}{\sqrt{d_k}} e_{\chi_0} & \mbox{if $\ff = \fq^l$}
\end{array}
\right.
$$
where $L^{(2)}(0, \tchi) := \lim\limits_{s \to 0} s^{-2} L(s,
\tchi)$. (Note in particular
that the Gauss Sums $g_{\fm(\chi)}(\tchi)$
disappear. For more details, see Lemma~5.1 of~\cite{twizas}.)
The values of $L^{(2)}(0, \tchi)$ can then be computed using
the method of \cite{DT}.

As for the $p$-adic computations,
since $k$ is quadratic, the results of
Section~\ref{quadformula} can be used to calculate
$\Phi_{\ff,T_p,p}(1)$
for any prime $p$ such that
\emph{$p$ is prime to $\ff$} and \emph{$p$ splits in $k$}.
Suppose also that \emph{$f$ divides $p-1$}.
This means that the additive character $j\circ\xi$ takes
values in $\mu_{p-1}$ hence in $\bbZp^\times$ for any $(\xi,I)\in W_\ff$.
Consequently, Corollary~\ref{cor:3B} implies that $\PhifTpp(1)\in\bbZ_p G$.
Moreover, the coefficients of the formal power series
$F(\uX;\xi,I,\rho,\varepsilon\rho)^j$
\etc\ lie in $\bbZp$ (since $I$ will be prime to $p$ and $p$ splits).
The assumption $f|(p-1)$ therefore speeds up the calculations considerably,
although it is unnecessary from the theoretical viewpoint and
places a major restriction on $p$. We shall assume from now
on that $p$ satisfies \emph{the three conditions above}.

The remarks at the beginning of Section~\ref{quadformula}
show that Corollary~\ref{cor:3B} now suffices in principle
for the numerical calculation
$\Phi_{\ff,T_p,p}^{(j)}(1)$. In practice, however,
the computation of the formal power series
$F(\uX;\xi,I,\rho,\varepsilon\rho)$ by means of\refeq{3B.5} can still
be prohibitively lengthy.  This is because the number of points $\ta$ in
$I\cap P(\rho,\varepsilon\rho)$ equals the index
$|I:\bbZ\rho+\bbZ\varepsilon\rho|$ which in turn is proportional to
the coefficients of $\varepsilon$ in a $\bbZ$-base $\{1,b\}$ of $\cO$
(for fixed, optimal $\rho$ and $I$).  But
these coefficients can be extremely large, even for $k$ of moderate discriminant and
(especially) $\ff$ of moderate norm. (Recall that $E_{\ff+}=\langle \varepsilon\rangle$.)
To tackle this problem the
approach of Corollary~\ref{cor:3B} can be refined as follows. Suppose
that $\rho_0,\ldots,\rho_L$ lie in $I\cap\kps$ but not in $\ker \xi$,
with $\rho_0=\rho$, $\rho_L=\varepsilon\rho$ and
\beql{eq:blob}
\mbox{
$\sgn\left(\det\left(
\begin{array}{c}
\uiota(\rho_{t-1})\\
\uiota(\rho_t)
\end{array}
\right)\right)=
\sgn\left(\det\left(
\begin{array}{c}
\uiota(\rho)\\
\uiota(\varepsilon\rho)
\end{array}
\right)\right)$
($=\sgn(\iota_2(\varepsilon)-\iota_1(\varepsilon))$)
for $t=1,\ldots,L$.}
\eeq
This condition means that the cone on $\rho$ and $\varepsilon\rho$ is
a `fan' of the cones on successive pairs $(\rho_{t-1},\rho_t)$.  More
precisely, $C(\rho,\varepsilon\rho)$ is the disjoint union
$C(\rho_{t-1},\rho_t)$ for $t=1,\ldots,L$ and it follows
from\refeq{3F.5} that $z_{T_p}(s;\xi,I,\rho,\varepsilon\rho)$ is the
sum of the $z_{T_p}(s;\xi,I,\rho_{t-1},\rho_t)$ for $\Re(s)>1$, hence
for all $s\in\bbC$ by analytic continuation. Thus the uniqueness of
the interpolation in Theorem~\ref{thm:3B}~\ref{part:thm3B1} implies
that
\beql{eq:3K}
z^{(j)}_{T_p,p}(s;\xi,I,\rho,\varepsilon\rho)=
\sum_{t=1}^L z^{(j)}_{T_p,p}(s;\xi,I,\rho_{t-1},\rho_t)
\eeq
for all $s\in D(p)$, and in particular for $s=1$.  (In fact,
the power series $F(\uX;\xi,I,\rho,\varepsilon\rho)$ equals $\sum_{t=1}^L
F(\uX;\xi,I,\rho_{t-1},\rho_t)$ but we do not need to know this.)  We
can therefore calculate $Z^{(j)}_{T_p,p}(1;\fw)$ by means of
Equations\refeq{3I.5}, and\refeq{3K}, using\refeq{3I.2} to determine
$z^{(j)}_{T_p,p}(1;\xi,I,\rho_{t-1},\rho_t)$ for each $t=1,2,\ldots,L$,
once the corresponding $F(\uX;\xi,I,\rho_{t-1},\rho_t)$ has been calculated.

Following work of Zagier (\cite{Zag}) Stark and others in similar contexts,
we now explain briefly how \emph{continued fractions}
may be used to obtain a sequence $\{\rho_t\}_{t=0}^L$ such that the index
$|I:\bbZ\rho_{t-1}+\bbZ\rho_t|$ is small for all $t$.
Most of the details can also be found in~\cite{Hay}.
Without loss of generality
we can assume that $\iota_1(\varepsilon)<1<\iota_2(\varepsilon)$.
Consider the following conditions on a pair of points $(x,y)\in (I\cap\kps)^2$
\[
\mbox{
(a)\ \
$\bbZ x+\bbZ y=I$,\ \ \ \
(b)\ \
$\iota_1(x)>\iota_1(y)$\ \ \ \
and (c)\ \
$\det\left(
\begin{array}{c}
\uiota(x)\\
\uiota(y)
\end{array}
\right)>0$
}
\]
It is easy to find a pair $(x,y)=(\trho_0,\trho_1)$, say, satisfying these conditions
and to see that they must then also hold with $(x,y)=(\trho_1,\trho_2)$ where
\[
\trho_2:=-\trho_0+b_1\trho_1
\ \ \mbox{and}\ \
b_1:=\lceil\iota_1(\trho_0/\trho_1)\rceil\geq 2
\]
hence also for $(x,y)=(\trho_2,\trho_3)$ where
$\trho_3:=-\trho_1+b_2\trho_2:=
-\trho_1+\lceil\iota_1(\trho_1/\trho_2)\rceil\trho_2$
and so on inductively.
In this way we produce an infinite sequence $\trho_0,\trho_1,\trho_2,\trho_3,\trho_4,\ldots$
such that Conditions (a)---(c) are obeyed for each successive pair
$(x,y)=(\trho_{n-1},\trho_n)$ for $n=1,2,3,\ldots$, and
$\trho_{n+1}=-\trho_{n-1}+b_n\trho_n$ where
$b_n:=\lceil\iota_1(\trho_{n-1}/\trho_n)\rceil\geq 2$.
In fact, we have a \emph{`type II continued fraction'} expansion
converging to $\iota_1(\trho_0/\trho_1)$:
\[
\iota_1(\trho_0/\trho_1)=
b_1-\frac{1}{{\displaystyle b_2}-{\displaystyle\frac{1}
{{\displaystyle b_3}-\ldots}}}
\]
The discreteness of $\uiota(I)$ implies that one cannot have
both $\iota_1(\trho_{n-1})>\iota_1(\trho_n)>0$
and also $\iota_2(\trho_{n-1})>\iota_2(\trho_n)>0$ indefinitely.
Thus there exists $N\geq 1$ such that $\iota_2(\trho_{N-1})<\iota_2(\trho_N)$
and one can show inductively that this property too is inherited from then on:
for each $n\geq N$ the pair
$(x,y)=(\trho_{n-1},\trho_n)$ must satisfy Condition~(a) together with
the following strengthening of Conditions~(b) and~(c)
\[
\mbox{
(b$'$)\ \
$\iota_1(x)>\iota_1(y)$\ \
and\ \ $\iota_2(x)<\iota_2(y)$}
\]
In the terminology of~\cite{Hay}, $\trho_{n-1}$ and $\trho_n$ are
successive points of the
\emph{`convexity polygon of $I$'}. But the group $E_+(k)$, hence
also $E_{\ff+}$, acts on this polygon and it follows that there exists $M>0$
such that
$\trho_{n+M}=\varepsilon\trho_n$ for all $n\geq N$. (This reflects the fact
that the sequence $\{b_n\}_{n=1}^\infty$ is eventually periodic.)
Thus, choosing any $n_0\geq N$ we obtain a finite
sequence $\{\rho'_m:=\trho_{n_0+m}\}_{m=0}^M$
with $\rho'_M=\varepsilon\rho'_0$ and
$|I:\bbZ\rho'_{m-1}+\bbZ\rho'_m|=1$ for $m=1,\ldots,M$.
Unfortunately, we may have $\rho'_m\in\ker(\xi)$ for some values of $m$,
but, since $\xi$ is non-trivial on $I$, such `bad' terms must
at least be \emph{non-consecutive}. We can therefore choose $n_0$
such that $\rho'_0\nin\ker(\xi)$ (so $\rho'_M\nin\ker(\xi)$)
and by simply skipping the bad terms and renumbering,
we finally arrive at a subsequence
$\rho_0=\rho'_0,\ldots,\rho_t,\ldots,\rho_L=\rho'_M=\varepsilon\rho_0$
in $(I\cap\kps)\setminus\ker(\xi)$ satisfying~\refeq{blob}.
Moreover, the indices $|I:\bbZ\rho_{t-1}+\bbZ\rho_t|=|I\cap P(\rho_{t-1},\rho_t)|$
are equal either to $1$ or to some
\emph{`partial quotient'} $b_n$ of the continued fraction.
Hence (empirically at least) they are still relatively small.
In fact one easily checks the explicit formula:
\[
I\cap P(\rho_{t-1},\rho_t)=
\left\{
\begin{array}{ll}
\{\rho'_{m-1}\}&
  \parbox{12 em}{if $(\rho_{t-1},\rho_t)=(\rho'_{m-1},\rho'_m)$
  for some $m$, $1\leq m\leq M$}\\
  &\\
\{\rho'_{m-1}\}\stackrel{\displaystyle\cdot}{\cup}
\{\rho'_m,2\rho'_m,\ldots,(b_{n_0+m}-1)\rho'_m\}&
  \parbox{12 em}{if $(\rho_{t-1},\rho_t)=(\rho'_{m-1},\rho'_{m+1})$
  for some $m$, $2\leq m\leq M-1$}\\
\end{array}
\right.
\]
which serves to calculate each $F(\uX;\xi,I,\rho_{t-1},\rho_t)$ and hence
$Z_{T_p,p}^{(j)}(1;\fw)$, as explained above.
In practice this leads to massive time-savings compared
with the use of $\rho_0$ and $\rho_L=\varepsilon \rho_0$ alone:
The effect of the smaller indices greatly outweighs
that of having $L$ such calculations instead of one.

We need to know how to ensure the accuracy of our calculated
value of $\PhifTpp(1)$ to a given number $N\geq 0$, say, of $p$-adic places.
By Equations\refeq{3K} and\refeq{3I.5}, it suffices to calculate
each value
$z^{(j)}_{T_p,p}(1;\xi,I,\rho_t,\rho_{t-1})$, $t=0,\ldots,L$
with an error less than $p^{-N}$ in $p$-adic absolute value.
For this, we fix $t$ and write $(1+\uX)\inv F(\uX;\xi,I,\rho_t,\rho_{t-1})$
as $\sum_{i,l\geq 0}a_{i,l}X_1^iX_2^l$. Consider
the real function $f_p(x)$ defined for all $x>-2$ by
$
f_p(x):=\frac{x+2}{p-1}-\frac{2}{\log p}\log(x/2+1)-2
$.
We note that $f_p$ is monotonic increasing and unbounded on the interval
$\Big[\frac{2(p-1)}{\log p}-2,\infty\Big)$. The following result
therefore solves our error-control problems.
\begin{prop}\label{prec}
With the above notations, choose
$M\geq \frac{2(p-1)}{\log p}-2$ such that $f_p(M)>N$ and suppose that
for each pair $(i,l)$ with $i,l\geq 0$ and $i+l<M$ we have computed
an element $\ta_{i,l}$ of $\bbZp[\mu_f]$ ($=\bbZp$)
such that $|\ta_{i,l}-a_{i,l}|_p<p^{-N}$. Then
\beql{eq:3S}
\left|z^{(j)}_{T_p,p}(1;\xi,I,\rho_t,\rho_{t-1})-
\frac{1}{p^2}\sum_{0\leq i,l\atop i+l<M}
\frac{c_{i+1}c_{l+1}\ta_{i,l}}{(i+1)(l+1)}\right|_p<p^{-N}
\eeq
\end{prop}
\bPf Lemma~\ref{lemma:estimcn} shows that
$c_{i+1}c_{l+1}/p^2(i+1)(l+1)$ is $p$-integral for
all $i,l\geq 0$. It is therefore
enough to show that Equation\refeq{3S} holds
with $\ta_{i,l}$ replaced by $a_{i,l}$ and by\refeq{3I.2} it suffices
to prove that $|c_{i+1}c_{l+1}a_{i,l}/p^2(i+1)(l+1)|_p<p^{-N}$
for any $i,l\geq 0$ with $i+l\geq M$. But the $a_{i,l}$ are
$p$-integral so the estimate of Lemma~\ref{lemma:estimcn}
together with the obvious estimate
$|1/(i+1)(l+1)|_p\leq (i+1)(l+1)\leq ((i+l+2)/2)^2$
shows that, for such $i,l$ we have
$|c_{i+1}c_{l+1}a_{i,l}/p^2(i+1)(l+1)|_p
\leq p^{-f_p(i+l)}\leq p^{-f_p(M)}<p^{-N}$ as required.
\ePf

In performing the calculations to compute (the approximations
$\ta_{i,j}$ to) the $a_{i,j}$, it is a good idea to
represent all power series
$A(\uX)\in\bbZp[[\uX]]$
as the sum of their homogeneous components $A_\nu(\uX)$ with $\nu
\geq 0$. Indeed, if we
are only interested in the coefficients of $X_1^iX_2^l$ with $i + l< M$ then
we can simply ignore the components
$A_\nu(\uX)$ with $\nu \geq M$. Moreover, if
$B(\uX)$ is another power series similarly represented as
$\sum_\nu B_\nu(\uX)$, then
the homogeneous components of the sum, product and
(assuming $p\ndiv B_0(\uX)$) the quotient of $A(\uX)$ by $B(\uX)$
can be calculated easily in terms of the $A_\nu(\uX)$
and $B_\nu(\uX)$. For instance, in the last case, if
$B_0(\uX)=b\in\bbZp^\times$
then $(A/B)_0(\uX)=b\inv A_0(\uX)$ and for $\nu\geq 1$ there is the simple
recurrence $(A/B)_\nu(\uX)=
b\inv[A_\nu(\uX)-B_1(\uX)(A/B)_{\nu-1}(\uX)-
\ldots-B_\nu(\uX)(A/B)_{0}(\uX))]$.

We next explain the basis of our method for working in $\WedtQUS$.
Let $\mathbf{X}=\{X_1, \dots, X_m\}$, say,
be the set of $\Gal(\bar\bbQ/\bbQ)$-conjugacy classes of characters in
$G^\ast$.
Then $\mathbf{X}$ can also be identified with the set of
isomorphism classes of irreducible rational representations of $G$,
a given conjugacy class $X_i\in\mathbf{X}$ corresponding
to the unique isomorphism class of representations
with character $\sum_{\chi\in X_i}\chi$. We write
$e_i$ for the rational idempotent $\sum_{\chi\in X_i}\chi\in \bbQ G$
of this character and $\bbQ(X_i)$ for $e_i\bbQ G$. Considered as a
$\bbQ G$-module, the latter is a representation lying in $X_i$.
Considered as a ring, $\bbQ(X_i)$ is a field and $\bbQ G$ is the direct product
$\prod_{i=1}^m\bbQ(X_i)$.
Thus we obtain a decomposition as $\bbQ G$-module
\beql{eq:dec}
\bbQ U_S=\bigoplus_{i=1}^m e_i\bbQ U_S\cong
\bigoplus_{i=1}^m\bbQ(X_i)^{r_i}
\eeq
where $r_i$ denotes the common value of $r(S, \chi)$ for all $\chi\in X_i$.
Clearly, for each $i$,
there exist $\bbQ(X_i)$-bases of $e_i\bbQ U_S$ consisting of
$S$-units (more precisely, elements of $1\otimes U_S$).
For small examples at least
it is not hard to find such a basis
by using the idempotents $e_i$. Let $v_{i,1},\ldots,v_{i,{r_i}}$
be such a basis, we then say that the elements
$v_{i,1},\ldots,v_{i,{r_i}}$ \emph{realize} the decomposition
(\ref{eq:dec}).
Passing to the exterior square, the product $v_{i,j}\wedge v_{i',j'}$
is zero in $\WedtQUS$ unless $i=i'$, so we obtain
a decomposition as $\bbQ G$-module
\[
\WedtQUS=\bigoplus_{i=1}^m e_i\WedtQUS=
\bigoplus_{i=1}^m\bigoplus_{1\leq j<j'\leq r_i}\bbQ(X_i)
(v_{i,j}\wedge v_{i,j'})
\]
Let $d_i$ denote the common order of each character $\chi\in X_i$.
The linear extension of any character $\chi\in X_i$ defines an
isomorphism $\chi:\bbQ(X_i)\rightarrow\bbQ(\mu_{d_i})$ of fields and of $\bbQ G$-modules
($G$ acting on $\bbQ(\mu_{d_i})$ \emph{via}
$\chi$). It follows that $\bbQ(X_i)$ has a $\bbQ$-basis of form
$\{1,\sigma_i\,\ldots,\sigma_i^{\phi(d_i)-1}\}$ where $\sigma_i\in G$
is any chosen element of $G$ such that $\chi(\sigma_i)$ is a primitive $d_i$th
root of unity. We thus obtain
$\bbQ$-bases $\cE$ of $\bbQ U_S$ and $\cB$ of $\WedtQUS$, of the forms
$\{\sigma_i^k v_{i,j}\}$ and $\{\sigma_i^k v_{i,j}\wedge v_{i,j'}\}$
respectively, where $i,j,j',k$ satisfy $1\leq i\leq m$,
$1\leq j<j'\leq r_i$ and $0\leq k\leq \phi(d_i)-1$.
A base for the subspace $\WedtQUS\St$ consists
of the subset $\cB_2=\{\sigma_i^k v_{i,1}\wedge v_{i,2}\,:\,r_i=2\}$ of $\cB$.
Having fixed these bases, we now explain how to find a
$\bbZ$-basis for the lattice
$\overline{\WedtZUS}\St$, expressed
as column vectors in the basis $\cB_2$.
First, let $\{u_1, \dots, u_t\}$ be any $\bbZ$-basis of the lattice
$1\otimes U_S\cong U_S/\{\pm 1\}$ in $\bbQ U_S$. We easily express
each $u_i$ as a rational linear combination of the base $\cE$.
By distributivity, each product $u_l\wedge u_{l'}$ is then expressed in the base $\cB$.
(For $k+k'\geq \phi(d_i)$, the product $\sigma_i^{k+k'} v_{i,j}\wedge v_{i,j'}$
can be reexpressed in the base $\cB$ by using the relations
$P(\sigma_i)v_{i,j}\wedge v_{i,j'}=0$
for any multiple $P$ in $\bbZ[X]$ of the $d_i$th cyclotomic polynomial.)
These products generate $\overline{\WedtZUS}$ over $\bbZ$. We form
the matrix $M$ of their rational column vectors in the base $\cB$,
with the coefficients of $\cB_2$ written first. Standard
column operations on $M$ reduce it to an Hermite Normal Form from
which we can read off a $\bbZ$-basis of the intersection
$\WedtQUS\St\cap\overline{\WedtZUS}=\overline{\WedtZUS}\St$
written in the $\bbQ$-basis $\cB_2$.
Exactly the same method can be used to find a $\bbZ$-basis for
the lattice $\overline{{\textstyle\bigwedge}^2_{\bbZ G}E(K)}\St$.
\subsection{The method of verification}\label{subsec:metover}
In this section, we explain the method used to numerically verify
Conjecture~\ref{conj:2F}. We illustrate this method using the first
example. Data on the verification of the conjecture in all the
examples (including the first one) are summarised in several tables
given at the end of this section. The first column of each table
contains the number of the example, the meaning of the other columns
of these tables is explained in the following subsections.

All the examples have been verified using the PARI/GP system
\cite{GP}.
\subsubsection{The extension $K/k$}
The data concerning the extension $K/k$ are summarised in
Table~1. First, we list the ground field $k$, its class number $h_k$,
the integral ideal $\ff$ and the conductor $\ff(K/k)$ dividing $\ff$
of the extension $K/k$.  (This is the minimal cycle modulo which $K$
is the ray-class field to $k$.)  The two integral ideals $\ff$ and
$\ff(K/k)$ are given as products of prime ideals in $k$, with $\fq_q$,
$\fq'_q$ denoting prime ideals in $k$ above the prime $q$ (if $q$ is
inert, we write $q\cO_k$ instead).  \textsl{In the first example, we
have $k := \bbQ(\sqrt{37})$, $\ff := 2\cO_k$ (\ie\ $2$ is inert in
$k$) and $\ff(K/k) = \ff$.}

Next, we give the monic irreducible polynomial $P_\theta(X)\in\bbZ[X]$
of an algebraic integer $\theta$ such that $K = \bbQ(\theta)$ (these
polynomials have been computed using the method of \cite{Ro}), the
factorisation of the discriminant of $K/\bbQ$, the class-number $h_K$
of $K$ and the structure of the Galois group $G$ as a product of
cyclic groups. Actually, in all examples but the last, the group $G$
is cyclic (it is isomorphic to $C_3 \times C_3$ in the last example)
and we let $\sigma$ denote a (fixed) generator of $G$ ($\sigma_1$ and
$\sigma_2$ are two (fixed) generators in the last example). The next
column of the table gives the action of this generator $\sigma$
(\emph{resp.}\ of $\sigma_1$ and $\sigma_2$) on the algebraic integer
$\theta$. In some examples, the expression for $\sigma(\theta)$ is too
long to be conveniently included in the table. Finally, the last entry
of the table is the degree $n_c$ of $K^c/K$ where $K^c$
is the Galois closure of $K/\bbQ$. Thus $K/\bbQ$ is Galois if and only
if $n_c = 1$.

\textsl{In the first example, we have $P_\theta(X) = X^6 - 3 X^5 - 2
X^4 + 9 X^3 - 5 X + 1$, $d_K = 2^4 \cdot 37^3$, $h_K = 1$, $G \simeq
C_3$, $\sigma(\theta) = -\theta^{5} + 2 \theta^{4} + 4 \theta^{3} - 6
\theta^{2} - 4 \theta + 3$ and $[K^c:K] = 1$.}
\subsubsection{The modules $\bbQ U_S$ and $\overline{\WedtZUS}\St$}
The corresponding data are summarised in Tables~2 and 3. We start
with the columns of Table~2. Whenever space allows it, the second column
gives a $\bbZ$-base of $1\otimes U_S\subset\bbQ U_S$.
We abuse notation by writing $u_i$ both for an
element of such a base and for the corresponding element of
$U_S$ itself (unique up to sign). The ranks of $1\otimes E(K)$ and
$1\otimes U_S\cong U_S/\{\pm 1\}$ are  $[K:\bbQ] - 1$ and
$[K:\bbQ] - 1 + |S_0|$ respectively. The first $[K:\bbQ] - 1$
elements of the given base lie in $1\otimes E(K)$.

\textsl{In the first example, the $\bbZ$-rank of $U_S$ is $6$ and the
rank of $E(K)$ is $5$. A $\bbZ$-base of $1\otimes E(K)$ is given by:
$u_{1} := \theta^{3} - 2 \theta^{2} - \theta + 1$, $u_{2} := \theta$,
$u_{3} := \theta^{5} - 2 \theta^{4} - 3 \theta^{3} + 4 \theta^{2} + 2
\theta - 1$, $u_{4} := \theta^{5} - 3 \theta^{4} - \theta^{3} + 7
\theta^{2} - 2 \theta - 1$ and $u_{5} := \theta^{5} - 2 \theta^{4} - 3
\theta^{3} + 5 \theta^{2} + 2 \theta - 2$. To get a base of $1\otimes
U_S$ we add to this system the element $u_{6} := \theta^{3} - 2
\theta^{2} - 2 \theta + 3$.}

Next, we give characters $\chi$ generating the group $G^\ast$ of
irreducible complex characters of $G$. These characters are
defined by their values on the specified generators of $G$.
Finally, we give the set
$\mathbf{X}=\{X_1, \dots, X_m\}$ considered as
irreducible rational characters of $G$.
Thus each $X_i$ is written as a sum of the elements of the corresponding
$\Gal(\bar\bbQ/\bbQ)$-conjugacy class in $G^\ast$.
In all the examples, the character $\chi_0$ denotes the
trivial character.

\textsl{In the first example, $G^\ast$ is generated by the character
$\chi$ with $\chi(\sigma) := e^{2i\pi/3}$ and there are two
irreducible rational characters $X_1 := \chi_0$ and $X_2 := \chi +
\chi^2$.}

We now look at the columns of Table~3. The second column contains the
structure of $\bbQ U_S$ as $\QG$-module as represented in
Equation~(\ref{eq:dec}). For the examples in which Table~2
lists the $\bbZ$-base of $1\otimes U_S$ used, the third column of Table~3
gives an isotypic $\bbQ$-basis
$\{v_{i,j}\}$ of $\bbQ U_S$, written relative to this base, which
`realizes' the decomposition in the first column (see Section~\ref{comp}).
Explicitly, the (integral) vector
$(a_1, \dots, a_l)$ represents the image of the
$S$-unit $\pm u_1^{a_1} \cdots u_l^{a_l}$ in $1\otimes U_S$.
The fourth column contains the
idempotent $\teSgt$. Note that $\overline{\WedtZUS} =
\overline{\WedtZUS}\St$ if and only if $\teSgt = 0$.

Finally, we give
an element $\gamma$ generating a $\ZG$-submodule of finite index
in $\overline{\WedtZUS}\St$, the index being given in the last
column. This element will be used below to verify the conjecture.
It can be found by looking among `small', random integral combinations
of the $\bbZ$-base of
$\overline{\WedtZUS}\St$ which in turn is found as
described in the last part of Section~\ref{comp}.
(Such a combination
generates a submodule of finite index if and only if,
when it is
expressed in the $\bbQ$-basis `$\cB_2$' of $\WedrQUS\St$ (see \emph{idem}),
the coefficient of $\sigma^k_i v_{i,1}\wedge v_{i,2}$ is non-zero
for some $k=k(i)$ for each $i$ with $r_i=2$.)
When the index is greater than $1$, it
is of course possible that there exists an element
$\gamma'$ generating a $\bbZ G$-submodule of smaller index in
$\overline{\WedtZUS}\St$.
We have expressed the element $\gamma$
as a sum of terms $(a_1, \dots, a_l) \wedge
(b_1, \dots, b_l)$ where each term $(a_1, \dots, a_l) \wedge (b_1,
\dots, b_l)$ represents the image of the element $u_1^{a_1} \cdots
u_l^{a_l} \wedge u_1^{b_1} \cdots u_l^{b_l}$ in
$\overline{\WedtZUS}$.

\textsl{In the first example, we have $\bbQ U_S \cong \bbQ(X_1)^2 +
\bbQ(X_2)^2$, isotypic elements realizing this decomposition are
$v_{1,1} := ( -2, 2, -1, -1, -3, 0)$, $v_{1,2} := ( -1, 2, -1, -2, 0,
3)$, $v_{2,1} := ( -2, 5, -4, -1, 0, 0)$ and $v_{2,2} := ( 3, -3, 3,
3, 0, 0)$. The idempotent $\teSgt$ vanishes so $\overline{\WedtZUS} =
\overline{\WedtZUS}\St$ and we find that $\gamma := ( 0, 0, 0, 0, 1,
0) \wedge ( 0, 0, 0, 0, 0, 1)$ generates $\overline{\WedtZUS}\St$ over
$\ZG$.}
\subsubsection{Verification of conjecture \ref{conj:2F}}
The data concerning the numerical verification of Conjecture~\ref{conj:2F}
are contained in Tables~4, 5, 6 and 7. More precisely, Tables~4 and
5 refer to parts (i), (iii) and the first statement of part~(iv) of
the conjecture while Tables~7 and 8 refer to part~(ii) and the second
statement of part~(iv).

Tables~4 and 5 give the same data for the first eight and the
last seven examples respectively. Their first two columns give computed
approximations to  $\frac{4}{\sqrt{d_k}}
R(\gamma)$ and to
$\Phi_{\ff,\emptyset}(1)$ (see the beginning of
Section~\ref{comp} for the computation of the latter. To save some
space, these values are given to a smaller precision than that to which
they were actually computed.) Now,
$\left(\WedtQUS\right)\St$ is free of rank 1 over $\bbQ G\St$,
generated by $\gamma$. It follows that if a solution $\eta_\ff$
of part~(i) of Conjecture~\ref{conj:2F} exists -- that is, if there
exists $\eta_\ff \in \left(\WedtQUS\right)\St$ such that
$
\frac{4}{\sqrt{d_k}} R(\eta_\ff) = \Phi_{\ff,\emptyset}(1)
$
-- then it must be of the form $\eta_\ff = A\gamma$ for some unique $A
\in \bbQ G\St$. So the equation to be solved becomes
\beql{eq:to_be_solved} A
\frac{4}{\sqrt{d_k}} R(\gamma) = \Phi_{\ff,\emptyset}(1)
\eeq
We solve this by first finding any solution $\tA$
of\refeq{to_be_solved} in $\bbR G$ to a high real precision.
This can be done using an obvious matrix method.
Applying $e_{S,2}$ to $\tA$ gives the solution in $\bbR G\St$ which
always turns out to be the approximation to the working precision of
an `obvious' element $A$ of $\bbQ G\St$.
The latter is listed in the fourth column
of Tables~4 and 5.

It is important to note that, whether or not the conjecture
holds, the non-vanishing of $R_K$ implies
that $\frac{4}{\sqrt{d_k}} R(\gamma)$ is \emph{invertible in $\bbR G\St$}
and hence that\refeq{to_be_solved} \emph{always has  a
unique solution in $\bbR G\St$}, namely
$(\sqrt{d_k}/4) R(\gamma)\inv \Phi_{\ff,\emptyset}(1)$.
On the other hand, parts~(iii) and~(iv) of the
conjecture predict that $A$ actually lies in $\frac{1}{b} \bbZ[1/g]G\St$
if $\ff \not= \fq^l$ (\emph{resp.}\ $\frac{1}{2b} \bbZ[1/g]G\St$ if $\ff =
\fq^l$) where $b$ denotes the index of $\bbZ G \gamma$
in $\overline{\WedtZUS}\St$. Even if the conjecture failed, such a solution could
always be `faked' to any desired real precision, simply
by approximating coefficients of $(\sqrt{d_k}/4) R(\gamma)\inv \Phi_{\ff,\emptyset}(1)$
sufficiently closely by elements of $\frac{1}{b}\bbZ[1/g]$ (or
$\frac{1}{2b} \bbZ[1/g]$).
Therefore, in order to verify the conjecture in a \emph{significant} way,
we must find that the `obvious' element $A\in\bbQ G\St$ determined above
has coefficients that lie naturally in $\frac{1}{b} \bbZ[1/g]G\St$
(or $\frac{1}{2b} \bbZ[1/g]G\St$). Moreover, their
\emph{numerators and denominators should be relatively
small and stabilise rapidly as the precision increases}.
This is indeed what we have observed in all our examples.

We set $\eta_\ff :=A\gamma$ as a
solution of part~(i) of the conjecture.
The next column of Tables~4 and 5 indicates whether the condition
$\ff = \fq^l$ applies in this example and the last column gives the
smallest positive integer $d_\ff$ such that $\eta_\ff$ belongs
$\frac{1}{d_\ff} \overline{\WedtZUS}\St$. This is determined
using the $\bbZ$-base of $\overline{\WedtZUS}\St$ found as described
in Section~\ref{comp}. According to part~(iii) of the
conjecture, we should have $d_\ff| g^e$ for some $e \geq 0$ if $\ff
\not= \fq^l$ (\ie\ if the answer in the previous column is `No'), and
by the first statement of part~(iv) of the conjecture $d_\ff|2g^e$
for some $e \geq 0$ if $\ff = \fq^l$ (\ie\ if the answer in the
previous column is `Yes'). Actually, in all examples, we have found
that $d_\ff = 1$ if $\ff \not= \fq^l$ and $d_\ff = 2$ if
$\ff =\fq^l$. (In the former case, therefore,
$\eta_\ff$ lies in $\overline{\WedtZUS}\St$. Indeed, using instead a
$\bbZ$-base of $\overline{\textstyle\bigwedge^2_{\bbZ G}E(K)}\St$
we have actually checked that it lies in this latter module,
see \cite[Rem. 3.4]{z1}.)

\textsl{We now illustrate this discussion using the first
example. We have computed
$$
\frac{4}{\sqrt{d_k}} R(\gamma)\! \simeq\!1.48595058394237662527436547684 +
    0.39280482164256390213294051602 (\sigma + \sigma^2)
$$
and
$$
\Phi_{\ff,\emptyset}(1) \simeq 0.35017047032862441050424222240 -
    0.74297529197118831263718273842 (\sigma + \sigma^2)
$$
A solution of\refeq{to_be_solved} in $\bbR G$ is then
$$
\tA\simeq 0.50000000000000000000000000000 -
0.50000000000000000000000000000 (\sigma + \sigma^2)
$$
Since $e_{S,2}$=1 in this example, we take $A$ to be the element
$\frac{1}{2}(1 - \sigma - \sigma^2)$ of $\frac{1}{2}\bbZ[1/3] G$ and
we set $\eta_\ff := A \gamma$. In this example, $\ff = \fq^l$ and we
find that $\eta_\ff$ belongs to $\frac{1}{2} \overline{\WedtZUS}\St$
but not $\overline{\WedtZUS}\St$ so $d_\ff = 2$. Hence, part~(i) and
the first statement of part~(iv) of Conjecture~\ref{conj:2F} are
numerically verified for this example up to the precision of the
computation.}\vertsp\\ \rem\ In this first example, as well as in
numbers 6, 7, 10, 13 and 15, it will be noticed that certain pairs of
coefficients coincide in $\Phi_{\ff,\emptyset}(1)$.  The explanation
is as follows.  Suppose for a moment that $k$ is any \emph{Galois}
extension of $\bbQ$ with $\Gamma=\Gal(k/\bbQ)$ and that the cycle
$\fm$ and the set $T$ are $\Gamma$-stable in the obvious sense. This
implies in particular that $k(\fm)$ is also Galois over $\bbQ$ and
$\Gamma$ acts by `extension and conjugation' on $G_\fm$ and $\bbR
G_\fm$. Explicitly, $\gamma(\sum_{g\in G_\fm}a_gg):= \sum_{g\in
G_\fm}a_g\tilde{\gamma}g\tilde{\gamma}\inv$ for any
$\tilde{\gamma}\in\Gal(k(\fm)/\bbQ)$ lifting any $\gamma\in\Gamma$. In
this situation, one can show that \beql{eq:xtra} \PhimT(s)=
\gamma(\PhimT(s))\ \ \ \forall\,\gamma\in\Gamma,\ \forall\,s\in\bbC,\
\Re(s)>1 \eeq By Theorem 2.3 of \cite{z1}, this equation can be
analytically continued to all $s\in\bbC\setminus\{1\}$ and even to
$s=1$ if Hypothesis~\ref{hyp:2A}~\ref{part:hyp2Apart3} holds.  All
these conditions are clearly met in the above-mentioned examples with
$\Gamma\cong\bbZ/2\bbZ$ and\refeq{xtra} therefore explains the
coincidence of coefficients in $\Phi_{\ff,\emptyset}(1)$, since in
each case $\Gamma$ acts by \emph{inversion} on $G=G_\ff$.  (A similar
coincidence in those of $\Phi_{\ff,T_p,p}(1)$ will follow by
interpolation from (the analytic continuation of) Equation\refeq{xtra}
for $s=m\in\cM(p)$.)  To prove Equation\refeq{xtra}, one uses obvious
actions of $\Gamma$ on $\fW_\fm$ and on $\Clmk$ in this situation (the
latter corresponds by the Artin map to the action on $G_\fm$) noting
that $\gamma(\fw_\fm^0)=\fw_\fm^0$ and that for all $\gamma\in\Gamma$,
$\fc\in\Clmk$, $\fw\in\fW_\fm$ and $\Re(s)>1$ we have
$\gamma(\fc\cdot\fw)=\gamma(\fc)\gamma(\fw)$ and, crucially,
$Z_T(s;\gamma(\fw))=Z_T(s;\fw)$. \vertsp\\ We now look at the entries
of Tables~6 and 7 which concern part~(ii) and the last statement of
part~(iv) of the conjecture. First, we list the prime numbers $p$ for
which part~(ii) of the conjecture has been tested. These prime numbers
$p$ must split in $k$ and satisfy $f|(p - 1)$ (which implies that
$(p,\ff)=1$).  Furthermore they must be relatively small for the
computation of $\Phi_{\ff, T_p, p}(1)$ to be feasible, \cf\
Proposition~\ref{prec}. Let $C_{p,\ff}$ be the constant appearing on
the L.H.S. of part~(ii) of the conjecture, \ie\ $C_{p,\ff} := 4 (1 -
p^{-1} \sigma_{\fp_1, \ff}) (1 - p^{-1} \sigma_{\fp_2, \ff})
j(\sqrt{d_k})^{-1}$ where $\fp_1$, $\fp_2$ are the two prime ideals in
$k$ above $p$. The next column contains the values of $C_{p,\ff}
R_p(\eta_\ff)$ for the $\eta_\ff$ computed above and for each of the
primes $p$. (The $p$-adic precision to which these are given is
smaller in most examples than the one that was used to verify the
conjecture.)  Each value of $C_{p,\ff} R_p(\eta_\ff)$ is checked
against the computed value of $\Phi_{\ff, T_p, p}(1)$. The two always
turn out to be equal, again up to the precision of the
computations. Thus, part~(ii) of the conjecture is satisfied to this
precision for these primes.

Note that $p$-adic numbers are written using the expansion to the
base $p$ with digits in the set $\{0,1,2,\ldots, p-1\}$.
The digits before the `decimal point'  correspond to negative powers
of $p$. If $p$ is larger than $10$ then we use the letters $A = 10, B
= 11, \dots$ to denote the extra digits. (The largest $p$ occurring is
$41$ for which we use the notation $(36), \dots, (40)$ to denote the
remaining digits.) The subscript at the end of the number is simply $p$.

\textsl{In the first example, part~(ii) of Conjecture~\ref{conj:2F}
has been numerically verified for $p = 3$, $7$ and $11$. As mentioned
above, we have found for each value of $p$ that $C_{p,\ff}
R_p(\eta_\ff) = \Phi_{\ff, T_p, p}(1)$ (up to our fixed $p$-adic
precision). We have found the following values}
\begin{center}
$\begin{array}{ll}
C_{3,\ff} R_{3}(\gamma) = &
   0.202021222001202022011122201001212121201_{3} \\
                          &
   {}+ 0.002112222212110120202010210110011000011_{3}
   (\sigma + \sigma^2)
\end{array}$

\vertsp

$C_{7,\ff} R_{7}(\gamma) = 0.232034003422155306164163_{7}
                 + 0.624214462041162660106331_{7} (\sigma + \sigma^2)$

\vertsp

$C_{11,\ff} R_{11}(\gamma) = 0.859AA8491A4592272_{11}
                 + 0.593A1A1A496337044_{11} (\sigma + \sigma^2)$
\end{center}

The next column contains the smallest positive integer $d_{\ff, \sigma
- 1}$ such that $(\sigma - 1)\eta_\ff$ belongs to $d_{\ff, \sigma -
1}^{-1} \overline{\WedtZUS}\St$ (in the last example, the smallest
positive integers $d_{\ff, \sigma_1 - 1}$ and $d_{\ff, \sigma_2 - 1}$
such that $(\sigma_i - 1)\eta_\ff$ belongs to $d_{\ff, \sigma_i -
1}^{-1} \overline{\WedtZUS}\St$ for $i = 1, 2$). Indeed, the second
statement of part~(iv) of the conjecture asserts that if $\ff = \fq^l$
then $I(\ZG) \eta_\ff \subset \bbZ[1/g] \overline{\WedtZUS}\St$. It
is easy to prove that if $G$ is generated by $\sigma$
(\emph{resp.}\ $\sigma_1$, $\sigma_2$ in the last example) then
$I(\ZG) = (\sigma - 1) \ZG$ (\emph{resp.}\ $I(\ZG) = (\sigma_1 - 1)\ZG +
(\sigma_2 - 1) \ZG$). Therefore, this statement is
true if and only if $(\sigma - 1) \eta_\ff$ (\emph{resp.}\ $(\sigma_1 - 1)
\eta_\ff$ and $(\sigma_2 - 1) \eta_\ff$) belongs to $\bbZ[1/g]
\overline{\WedtZUS}\St$, that is,
if $d_{\ff, \sigma - 1}$ (\emph{resp.}\ $d_{\ff,
\sigma_1 - 1}$ and $d_{\ff, \sigma_2 - 1}$) divide $g^e$ for
some $e \geq 0$. In fact, in all examples with $\ff = \fq^l$ we have
found that $d_{\ff, \sigma - 1}$ (\emph{resp.}\ $d_{\ff, \sigma_1 - 1}$ and
$d_{\ff, \sigma_2 - 1}$) is actually $1$. In other
words, Condition (49) of \cite{z1} is verified for these
examples. Indeed, we have actually checked that \cite[eq.\ (50)]{z1}
is verified (see Remark 3.4, \emph{ibid.}).

Finally, the last column gives the index of $\ZG \eta_\ff$ in
$d_\ff^{-1} \overline{\WedtZUS}\St$. In each case, we have found that
it is a small power of $2$ if $\ff = \fq^l$ and that $\eta_\ff$
actually generates $\overline{\WedtZUS}\St$ if $\ff \not= \fq^l$. We
do not expect this last fact to generalise.
There is no reason to expect $\overline{\WedtZUS}\St$
to be cyclic over $\bbZ G$ in general, and even when it is, the index of
$\bbZ G \eta_\ff$ should reflect the class number $h_K$ of $K$
as is the case with cyclotomic units, for $k=\bbQ$. Thus this index
might well be non-trivial in `larger' examples. Nevertheless, the mere
fact in all our examples this index is always very small (if not trivial)
\emph{is} significant: we certainly would not expect this if
$(\sqrt{d_k}/4) R(\gamma)\inv \Phi_{\ff,\emptyset}(1)$ were a random element
of $\bbR G\St$ and $A$ a `faked' approximation to the given precision,
lying in
$\frac{1}{b} \bbZ[1/g]G\St$ or $\frac{1}{2b} \bbZ[1/g]G\St$
(see discussion above).

\textsl{In the first example, $\ZG \eta_\ff$ is of index $4$ in
$\frac{1}{2} \overline{\WedtZUS}\St$}.

\begin{sidewaystable}
\centering
{\fontsize{8pt}{9pt}\selectfont
\begin{tabular}{|c|c|c|c|c|c|c|c|c|c|c|}
\hline
\onecolc{11}{\footnotesize \textrm{Table 1:} \emph{The extension $K/k$
\phv{12}{6}}} \\
\hline
 \# \phv{11}{5}
    & $k$   & $h_k$ & $\ff$ & $\ff(K/k)$ & $P_\theta(X)$ & $d_K$
    & $h_K$ & $G$   & $\sigma(\theta)$ & $n_c$ \\
\hline
1  \phv{11}{5}
   & $\bbQ(\sqrt{37})$
   & $1$
   & $2\cO_k$
   & $2\cO_k$
   & $X^{6} - 3 X^{5} - 2 X^{4} + 9 X^{3} - 5 X + 1$
   & $2^{4} \cdot 37^{3}$
   & $1$
   & $C_3$
   & $-\theta^{5} + 2 \theta^{4} + 4 \theta^{3} - 6 \theta^{2} - 4
      \theta + 3$
   & $1$    \\
\hline
2  \phv{11}{5}
   & $\bbQ(\sqrt{43})$
   & $1$
   & $\fq_3^2$
   & $\fq_3^2$
   & $X^{6} - 16 X^{4} - 12 X^{3} + 21 X^{2} + 10X - 7$
   & $2^{6} \cdot 3^{4} \cdot 43^{3}$
   & $1$
   & $C_3$
   & $\frac{1}{19}(7 \theta^{5} - 3 \theta^{4} - 108 \theta^{3} - 35 \theta^{2} + 105 \theta - 13)$
   & $3$    \\
\hline
3  \phv{11}{5}
   & $\bbQ(\sqrt{82})$
   & $4$
   & $2$
   & $\cO_k$
   & \onecoll{1}{$X^{8} + 2 X^{7} - 21 X^{6} - 78 X^{5} - 53 X^{4}$}
   & $2^{12} \cdot 41^{4}$
   & $1$
   & $C_4$
   & not given (too big)
   & $1$    \\
   \phv{5}{5}
   & & & & 
   & \onecolr{1}{${}+ 88 X^{3} + 114 X^{2} + 24 X - 4$}
   & & & & & \\
\hline
4  \phv{11}{5}
   & $\bbQ(\sqrt{89})$
   & $1$
   & $\fq_5$
   & $\fq_5$
   & $X^{4} + 2 X^{3} - 8 X^{2} - 9 X - 2$
   & $5 \cdot 89^{2}$
   & $1$
   & $C_2$
   & $-\theta - 1$
   & $2$    \\
\hline
5  \phv{11}{5}
   & $\bbQ(\sqrt{321})$
   & $3$
   & $\fq_2$
   & $\cO_k$
   & $X^{6} + 2 X^{5} - 18 X^{4} - 55 X^{3} - 26 X^{2} + 21 X + 3$
   & $3^{3} \cdot 107^{3}$
   & $1$
   & $C_3$
   & $\frac{1}{18}(- \theta^{5} + 3 \theta^{4} + 15 \theta^{3} - 26
      \theta^{2} - 48 \theta + 3)$
   & $1$    \\
\hline
6  \phv{11}{5}
   & $\bbQ(\sqrt{349})$
   & $1$
   & $2\cO_k$
   & $2\cO_k$
   & $X^{6} + 3 X^{5} - 36 X^{4} - 77 X^{3} + 200 X^{2} + 239 X - 205$
   & $2^{4} \cdot 349^{3}$
   & $1$
   & $C_3$
   & \onecoll{1}{$\frac{1}{385}(5 \theta^{5} + 2 \theta^{4} - 216
                  \theta^{3} + 38 \theta^{2}$}
   & $1$    \\
   \phv{5}{5}
   & & & & & & & &
   & \onecolr{1}{${}+ 1856 \theta - 335)$}
   &        \\
\hline
7  \phv{11}{5}
   & $\bbQ(\sqrt{401})$
   & $5$
   & $\fq_2\fq_2'$
   & $\cO_k$
   & \onecoll{1}{$X^{10} + 2 X^{9} - 20 X^{8} - 2 X^{7} + 69 X^{6} +
                  X^{5}$}
   & $401^{5}$
   & $1$
   & $C_5$
   & \onecoll{1}{$\frac{1}{27}(-7 \theta^{9} - 8 \theta^{8} + 151
                  \theta^{7} - 106 \theta^{6} - 473 \theta^{5}$}
   & $1$    \\
   \phv{5}{5}
   & & & &
   & \onecolr{1}{${}- 69 X^{4} - 2 X^{3} + 20 X^{2} + 2 X - 1$}
   & & &
   & \onecolr{1}{${}+ 359 \theta^{4} + 427 \theta^{3} - 220 \theta^{2}
                    - 67 \theta + 7)$}
   &        \\
\hline
8  \phv{11}{5}
   & $\bbQ(\sqrt{401})$
   & $5$
   & $\fq_5$
   & $\fq_5$
   & \onecoll{1}{$X^{20} + 2 X^{19} - 27 X^{18} - 58 X^{17} + 272
X^{16} + 639 X^{15}$}
   & $5^{5} \cdot 401^{10}$
   & $1$
   & $C_{10}$
   & not given (too big)
   & $2$    \\
   \phv{5}{5}
   & & & &
   & \onecol{1}{${}- 1245 X^{14} - 3339 X^{13} + 2469 X^{12}+ 8464 X^{11}$}
   & & & & & \\
   \phv{5}{5}
   & & & &
   & \onecol{1}{\quad ${} - 1650 X^{10} - 9965 X^{9} + 827 X^{8} + 6081 X^{7}$}
   & & & & & \\
   \phv{5}{5}
   & & & &
   & \onecolr{1}{${}- 914 X^{6} - 1796 X^{5} + 510 X^{4} + 151 X^{3} -
63 X^{2} + X + 1$} 
   & & & & & \\
\hline
9  \phv{11}{5}
   & $\bbQ(\sqrt{577})$
   & $7$
   & $\fq_2$
   & $\cO_k$
   & \onecoll{1}{$X^{14} + 2 X^{13} - 25 X^{12} - 69 X^{11} + 161
X^{10} + 632 X^{9}$} 
   & $577^{7}$
   & $1$
   & $C_7$
   & not given (too big)
   & $1$    \\
   \phv{5}{5}
   & & & &
   & \onecol{1}{${}- 147 X^{8} - 2146 X^{7} - 1171 X^{6} + 2669 X^{5}
+ 2682 X^{4}$} 
   & & & & & \\
   \phv{5}{5}
   & & & &
   & \onecolr{1}{${} -  667 X^{3} - 1466 X^{2} - 336 X + 49$}
   & & & & & \\
\hline
10 \phv{11}{5}
   & $\bbQ(\sqrt{709})$
   & $1$
   & $2\cO_k$
   & $2\cO_k$
   & $X^{6} - 56 X^{4} + 784 X^{2} - 2836$
   & $2^{4} \cdot 709^{3}$
   & $1$
   & $C_3$
   & $\frac{1}{212}(-9 \theta^{4} + 420 \theta^{2} - 106 \theta -
      3136)$
   & $1$    \\
\hline
11 \phv{11}{5}
   & $\bbQ(\sqrt{709})$
   & $1$
   & $2\fq_5$
   & $2\fq_5$
   & \onecoll{1}{$X^{12} - 53 X^{10} + 970 X^{8} - 7657 X^{6}$}
   & $2^{8} \cdot 5^{3} \cdot 709^{6}$
   & $1$
   & $C_6$
   & \onecoll{1}{$\frac{1}{405650}(-231 \theta^{11} + 10953 \theta^{9}
                  - 164825\theta^{7}$}
   & $2$    \\
   \phv{5}{5}
   & & & &
   & \onecolr{1}{${}+ 25350 X^{4}- 29025 X^{2} + 6125$}
   & & &
   & \onecolr{1}{${}+ 920367 \theta^{5} - 1325445 \theta^{3} -
                  733225\theta)$}
   & \\
\hline
12 \phv{11}{5}
   & $\bbQ(\sqrt{1021})$
   & $1$
   & $\fq_5$
   & $\fq_5$
   & $X^{4} + 2 X^{3} - 32 X^{2} - 33 X + 17$
   & $5 \cdot 1021^{2}$
   & $1$
   & $C_2$
   & $-\theta - 1$
   & $2$    \\
\hline
13 \phv{11}{5}
   & $\bbQ(\sqrt{2069})$
   & $1$
   & $2\cO_k$
   & $2\cO_k$
   & $X^{6} - 84 X^{4} + 1764 X^{2} - 8276$
   & $2^{4} \cdot 2069^{3}$
   & $1$
   & $C_3$
   & $\frac{1}{60}(\theta^{4} - 70 \theta^{2} - 30 \theta + 784)$
   & $1$    \\
\hline
14 \phv{11}{5}
   & $\bbQ(\sqrt{2069})$
   & $1$
   & $2\fq_5$
   & $2\fq_5$
   & \onecoll{1}{$X^{12} - 71 X^{10} - 134 X^{9} + 1128 X^{8}  + 3138
X^{7}$} 
   & $2^{8} \cdot 5^{3} \cdot 2069^{6}$
   & $1$
   & $C_6$
   & not given (too big)
   & $2$    \\
   \phv{5}{5}
   & & & &
   & \onecol{1}{${}- 2847 X^{6} - 12804 X^{5} - 2686 X^{4} + 13110 X^{3}$} 
   & & & & & \\
   \phv{5}{5}
   & & & &
   & \onecolr{1}{${}+ 9935 X^{2} + 2150 X + 125$}
   & & & & & \\
\hline
15 \phv{11}{5}
   & $\bbQ(\sqrt{9897})$
   & $3$
   & $\fq_3^2$
   & $\fq_3^2$
   & \onecoll{1}{$X^{18} - 204 X^{16} + 15822 X^{14} - 590238 X^{12}$} 
   & $3^{21} \cdot 3299^{9}$
   & $3$
   & $C_3^2$
   & \onecoll{1}{$G$ has two generators $\sigma_1$ and $\sigma_2$
     which}
   & $1$    \\
   \phv{5}{5}
   & & & &
   & \onecol{1}{${}+ 11246949 X^{10} - 110721114 X^{8} + 550866177 X^{6}$}
   & & &
   & \onecoll{1}{are not given (too big)}
   & \\
   \phv{5}{5}
   & & & &
   & \onecolr{1}{${}- 1324310688 X^{4} + 1327290624 X^{2} - 364843008$}
   & & & & & \\
\hline
\end{tabular}
}
\end{sidewaystable}

\begin{sidewaystable}
\centering
{\fontsize{8pt}{9pt}\selectfont
\begin{tabular}{|c|c|c|c|}
\hline
\onecolc{4}{\footnotesize \textrm{Table 2:} \emph{The modules
$\WedtZUS$ and $\overline{\WedtZUS}\St$} \phv{18}{6}} \\
\hline
 \# \phv{11}{5}
    & $S$-units & $G^\ast$ & $\mathbf{X}$ \\
\hline
1  \phv{11}{5}
   & $u_{1} := \theta^{3} - 2 \theta^{2} - \theta + 1, \ \
      u_{2} := \theta$
   & $\langle \chi \rangle$ with $\chi(\sigma) := e^{2i\pi/3}$
   & $X_1 := \chi_0$ \\
   \phv{5}{5}
   & $u_{3} := \theta^{5} - 2 \theta^{4} - 3 \theta^{3} + 4 \theta^{2}
               + 2 \theta - 1,\ \
      u_{4} := \theta^{5} - 3 \theta^{4} - \theta^{3} + 7 \theta^{2} -
               2 \theta - 1$
   &
   & $X_2 := \chi + \chi^2$ \\
   \phv{5}{5}
   & $u_{5} := \theta^{5} - 2 \theta^{4} - 3 \theta^{3} + 5 \theta^{2}
               + 2 \theta - 2,\ \
      u_{6} := \theta^{3} - 2 \theta^{2} - 2 \theta + 3$
   & &      \\
\hline
2  \phv{11}{5}
   & not given (too big)
   & $\langle \chi \rangle$ with $\chi(\sigma) := e^{2i\pi/3}$
   & $X_1 := \chi_0$, $X_2 := \chi + \chi^2$ \\
\hline
3  \phv{11}{5}
   & not given (too big)
   & $\langle \chi \rangle$ with $\chi(\sigma) := i$
   & $X_1 := \chi_0$, $X_2 := \chi^2$, $X_3 = \chi + \chi^3$ \\
\hline
4  \phv{11}{5}
   & $u_{1} := \theta^{3} + \theta^{2} - 10 \theta - 3,\ \
      u_{2} := \theta^{3} + \theta^{2} - 8 \theta - 3$
   & $\langle \chi \rangle$ with $\chi(\sigma) := -1$
   & $X_1 := \chi_0$ \\
   \phv{5}{5}
   & $u_{3} := \theta^{3} + 6 \theta^{2} + 5 \theta + 1,\ \
      u_{4} := 2 \theta + 1$
   &
   & $X_2 := \chi$ \\
\hline
5  \phv{11}{5}
   & $u_{1} := \frac{1}{12}(\theta^{5} - \theta^{4} - 19 \theta^{3} -
               2 \theta^{2} + 48 \theta + 9)$,\ \ 
     $u_{2} := \frac{1}{36}(\theta^{5} + 3 \theta^{4} - 27 \theta^{3}
               - 58 \theta^{2} + 84 \theta - 3)$,
   & $\langle \chi \rangle$ with $\chi(\sigma) := e^{2i\pi/3}$
   & $X_1 := \chi_0$ \\
   \phv{5}{5}
   & $u_{3} := \frac{1}{18}(2 \theta^{5} + 3 \theta^{4} - 30
               \theta^{3} - 101 \theta^{2} - 111 \theta - 15)$,\ \ 
     $u_{4} := \frac{1}{2}(\theta^{4} - 17 \theta^{2} - 21 \theta -
               1)$,
   &
   & $X_2 := \chi + \chi^2$ \\
   \phv{5}{5}
   & $u_{5} := \frac{1}{18}(\theta^{5} - 21 \theta^{3} - 25 \theta^{2}
               + 3 \theta + 6)$,\ \ 
     $u_{6} := \frac{1}{4}(\theta^{5} + \theta^{4} - 19 \theta^{3} -
               36 \theta^{2} + 10 \theta + 11)$
   & &      \\
\hline
6  \phv{11}{5}
   & not given (too big)
   & $\langle \chi \rangle$ with $\chi(\sigma) := e^{2i\pi/3}$
   & $X_1 := \chi_0$,
     $X_2 := \chi + \chi^2$ \\
\hline
7  \phv{11}{5}
   & not given (too big)
   & $\langle \chi \rangle$ with $\chi(\sigma) := e^{2i\pi/5}$
   & $X_1 := \chi_0$ \\
   \phv{5}{5}
   &
   &
   & $X_2 := \chi + \chi^2 + \chi^3 + \chi^4$ \\
\hline
8  \phv{11}{5}
   & not given (too big)
   & $\langle \chi \rangle$ with $\chi(\sigma) := e^{2i\pi/10}$
   & $X_1 := \chi_0$, $X_2 := \chi^5$ \\
   \phv{5}{5}
   &
   &
   & $X_3 := \chi^2 + \chi^4 + \chi^6 + \chi^8$ \\
   \phv{5}{5}
   &
   &
   & $X_4 := \chi + \chi^3 + \chi^7 + \chi^9$ \\
\hline
9  \phv{11}{5}
   & not given (too big)
   & $\langle \chi \rangle$ with $\chi(\sigma) := e^{2i\pi/7}$
   & $X_1 := \chi_0$ \\
   \phv{5}{5}
   &
   &
   & $X_2 := \chi + \chi^2 + \chi^3 + \chi^4 + \chi^5 + \chi^6$ \\
\hline
10 \phv{11}{5}
   & not given (too big)
   & $\langle \chi \rangle$ with $\chi(\sigma) := e^{2i\pi/3}$
   & $X_1 := \chi_0$, $X_2 := \chi + \chi^2$ \\
\hline
11 \phv{11}{5}
   & not given (too big)
   & $\langle \chi \rangle$ with $\chi(\sigma) := e^{2i\pi/6}$
   & $X_1 := \chi_0$, $X_2 := \chi^3$ \\
   \phv{5}{5}
   &
   &
   & $X_3 := \chi^2 + \chi^4$, $X_4 := \chi + \chi^5$ \\
\hline
12 \phv{11}{5}
   & $u_{1} := \frac{1}{9}(5 \theta^{3} + 3 \theta^{2} - 166 \theta +
               62),\ \
      u_{2} := \frac{1}{9}(4 \theta^{3} + 24 \theta^{2} - 14 \theta -
               53)$
   & $\langle \chi \rangle$ with $\chi(\sigma) := -1$
   & $X_1 := \chi_0$ \\
   \phv{5}{5}
   & $u_{3} := \frac{1}{3}(938 \theta^{3} - 4029 \theta^{2} - 5458
               \theta + 2600),\ \
      u_{4} := \frac{1}{9}(8 \theta^{3} + 57 \theta^{2} + 44 \theta -
               25)$
   &
   & $X_2 := \chi$ \\
\hline
13 \phv{11}{5}
   & not given (too big)
   & $\langle \chi \rangle$ with $\chi(\sigma) := e^{2i\pi/3}$
   & $X_1 := \chi_0$, $X_2 := \chi + \chi^2$ \\
\hline
14 \phv{11}{5}
   & not given (too big)
   & $\langle \chi \rangle$ with $\chi(\sigma) := e^{2i\pi/6}$
   & $X_1 := \chi_0$, $X_2 := \chi^3$ \\
   \phv{5}{5}
   &
   &
   & $X_3 := \chi^2 + \chi^4$, $X_4 := \chi + \chi^5$ \\
\hline
15 \phv{11}{5}
   & not given (too big)
   & \onecoll{1}{$\langle \chi_1, \chi_2 \rangle$ with
                 $\chi_1(\sigma_1) := e^{2i\pi/3}$,}
   & $X_1 := \chi_0$ \\
   \phv{5}{5}
   &
   & \onecol{1}{$\chi_1(\sigma_2) := 1$, $\chi_2(\sigma_1) := 1$}
   & $X_2 := \chi_1 + \chi_1^2$, $X_3 := \chi_2 + \chi_2^2$ \\
   \phv{5}{5}
   &
   & \onecolr{1}{and $\chi_2(\sigma_2) := e^{2i\pi/3}$}
   & $X_4 := \chi_1\chi_2 + \chi_1^2\chi_2^2$,
     $X_5 := \chi_1\chi_2^2 + \chi_1^2\chi_2$ \\
\hline
\end{tabular}
}
\end{sidewaystable}

\begin{sidewaystable}
\centering
{\fontsize{8pt}{10pt}\selectfont
\begin{tabular}{|c|c|c|c|c|c|}
\hline
\onecolc{6}{\footnotesize\textrm{Table 3:} \emph{The modules $\WedtZUS$
and $\overline{\WedtZUS}\St$ \phv{18}{6}}} \\
\hline
 \# \phv{12}{6}
    & $\bbQ U_S$ & generators of $\WedtQUS$
    & $\teSgt$
    & $\gamma$ & index of $\ZG\gamma$ \\
\hline
1  \phv{12}{6}
    & $\bbQ(X_1)^2 + \bbQ(X_2)^2$
    & $v_{1,1} := (-2,  2, -1, -1, -3, 0),$
    & $0$
    & $( 0,  0,  0,  0,  1,  0) \wedge  ( 0,  0,  0,  0,  0,  1)$
    & $1$           \\
    \phv{6}{6}
    &
    & $v_{1,2} := (-1,  2, -1, -2,  0, 3),$
    & & &       \\
    \phv{6}{6}
    &
    & $v_{2,1} := (-2,  5, -4, -1,  0,  0),$
    & & &       \\
    \phv{6}{6}
    &
    & $v_{2,2} := ( 3, -3,  3,  3,  0,  0)$
    & & &       \\
\hline
2  \phv{12}{6}
    & $\bbQ(X_1)^2  + \bbQ(X_2)^2$ 
    & not given
    & $0$
    & not given
    & $1$           \\
\hline
3  \phv{12}{6}
    & $\bbQ(X_1)^2  + \bbQ(X_2)^3 + \bbQ(X_3)^2$ 
    & not given
    & $1  - \sigma + \sigma^2 - \sigma^3$
    & not given
    & $1$           \\
\hline
4  \phv{12}{6}
   & $\bbQ(X_1)^2 + \bbQ(X_2)^2$
   & $v_{1,1} := ( 1,  0, -2,  0),\
      v_{1,2} := ( 0,  0,  0,  2)$
   & $0$
   & \onecoll{1}{$( 0,  1,  0,  0) \wedge  ( 0,  0,  -1,  0)$}
   & $2$        \\
   \phv{6}{6}
   &
   & $v_{2,1} := ( 1, -2,  0,  0),\
      v_{2,2} := ( 0, -2,  0,  0)$
   &
   & \onecolr{1}{${}+  ( 0,  0,  1,  0) \wedge ( 0,  0,  0,  -1)$}
   &            \\
\hline
5  \phv{12}{6}
   & $\bbQ(X_1)^2 + \bbQ(X_2)^2$
   & $v_{1,1} := ( 0,  0, -1, -1, -3,  0),$
   & $0$
   & $( 0,  0,  0,  0,  1,  0) \wedge  ( 0,  0,  0,  0,  0,  1)$
   & $1$        \\
   \phv{6}{6}
   &
   & $v_{1,2} := ( 1, -2,  0,  0,  0,  3),$
   & & &        \\
   \phv{6}{6}
   &
   & $v_{2,1} := ( 0,  3,  0,  0,  0,  0),$
   & & &        \\
   \phv{6}{6}
   &
   & $v_{2,2} := ( 0,  0,  3,  0,  0,  0)$
   & & &        \\
\hline
6  \phv{12}{6}
   & $\bbQ(X_1)^2 + \bbQ(X_2)^2$
   & not given
   & $0$
   & not given
   & $1$        \\
\hline
7  \phv{12}{6}
   & $\bbQ(X_1)^3 + \bbQ(X_2)^2$
   & not given
   & $1 + \sigma + \sigma^2 + \sigma^3 + \sigma^4$
   & not given
   & $1$        \\
\hline
8  \phv{12}{6}
   & $\bbQ(X_1)^2 + \bbQ(X_2)^2 + \bbQ(X_3)^2 + \bbQ(X_4)^2$
   & not given
   & $0$
   & not given
   & $2 \cdot 41$   \\
\hline
9  \phv{12}{6}
   & $\bbQ(X_1)^2 + \bbQ(X_2)^2$
   & not given
   & $0$
   & not given
   & $1$        \\
\hline
10 \phv{12}{6}
   & $\bbQ(X_1)^2 + \bbQ(X_2)^2$
   & not given
   & $0$
   & not given
   & $1$        \\
\hline
11 \phv{12}{6}
   & $\bbQ(X_1)^3 + \bbQ(X_2)^2 + \bbQ(X_3)^2 + \bbQ(X_4)^2$
   & not given
   & $1  + \sigma + \sigma^2 + \sigma^3 + \sigma^4 + \sigma^5$
   & not given
   & $1$        \\
\hline
12 \phv{12}{6}
   & $\bbQ(X_1)^2 + \bbQ(X_2)^2$
   & $v_{1,1} := (-1, -1, -2,  0),\
      v_{1,2} := (-1, -1,  0,  2)$
   & $0$
   & \onecoll{1}{$( 1,  0,  0,  0) \wedge  ( 0,  0,  -1,  0)$}
   & $2$        \\
   \phv{6}{6}
   &
   & $v_{2,1} := ( 1,  1,  0,  0),\
      v_{2,2} := ( 0,  2,  0,  0)$
   &
   & \onecolr{1}{${}+  ( 0,  0,  1,  0) \wedge ( 0,  0,  0,  1)$}
   &            \\
\hline
13 \phv{12}{6}
   & $\bbQ(X_1)^2 + \bbQ(X_2)^2$
   & not given
   & $0$
   & not given
   & $1$        \\
\hline
14 \phv{12}{6}
   & $\bbQ(X_1)^3 + \bbQ(X_2)^2 + \bbQ(X_3)^3 + \bbQ(X_4)^2$
   & not given
   & $3  + 3\sigma^3$
   & not given
   & $1$        \\
\hline
15 \phv{12}{6}
   & \onecoll{1}{$\bbQ(X_1)^2 + \bbQ(X_2)^2 + \bbQ(X_3)^3$~~}
   & not given
   & $(1 + \sigma_1 + \sigma_1^2)(2 - \sigma_2 - \sigma_2^2)$
   & not given
   & $1$        \\
\hline
\end{tabular}
}
\end{sidewaystable}

\begin{sidewaystable}
\centering
{\fontsize{8pt}{9pt}\selectfont
\begin{tabular}{|c|c|c|c|c|c|}
\hline
\onecolc{6}{\footnotesize \textrm{Table 4:} \emph{Verification of
Conjecture 2.2 - {\rm part (i), (iii) and (iv)}
\phv{11}{5}}} \\
\hline
\# \phv{11}{5}
   & \makebox[6.5cm][c]{$4\sqrt{d_k}\inv R(\gamma)$}
   & \makebox[6.5cm][c]{$\Phi_{\ff,\emptyset}(1)$}
   & \makebox[6cm][c]{$A$}
   & is $\ff = \fq^l$?
   & $d_\ff$        \\
\hline
1  \phv{11}{5}
   & \onecoll{1}{$1.4859505839423766252743654$}
   & \onecoll{1}{$0.3501704703286244105042422$}
   & $\frac{1}{2} (1  - \sigma - \sigma^2)$
   & Yes
   & $2$        \\
   \phv{5}{5}        
   & \onecolr{1}{${}+ 0.3928048216425639021329405 (\sigma +
                                                       \sigma^2)$}
   & \onecolr{1}{${}- 0.7429752919711883126371827 (\sigma +
                                                       \sigma^2)$}
   & & &        \\
\hline
2  \phv{11}{5}
   & \onecoll{1}{${}-  1.443448363709350198869637$}
   & \onecoll{1}{$     0.366471740478024869658988$}                    
   & $\frac{1}{2} (1  - \sigma + \sigma^2)$
   & Yes
   & $2$        \\
   \phv{5}{5}        
   & \onecol{1}{${}+  0.327469144766235117326033 \sigma$} 
   & \onecol{1}{${}-  0.039002595711789752332954 \sigma$} 
   & & &        \\
   \phv{5}{5}        
   & \onecolr{1}{${}-  1.848922699899164820861579 \sigma^2$}
   & \onecolr{1}{${}-  1.809920104187375068528625 \sigma^2$}
   & & &        \\
\hline
3  \phv{11}{5}
   & \onecoll{1}{$     0.990736966647953569158853 (1 + \sigma^3)$}
   & \onecoll{1}{$     0.326298344657731059802665$}
   & $\frac{1}{4} (-1  +  2\sigma -  3\sigma^2)$
   & Yes
   & $2$        \\
   \phv{5}{5}        
   & \onecolr{1}{${}-  0.104818803994323556683935 (\sigma + \sigma^2)$}
   &  \onecol{1}{${}-  0.221479540663407503118729 (\sigma + \sigma^3)$}
   & & &        \\
   \phv{5}{5}        
   & & 
     \onecolr{1}{${}-  0.769257425984546066040123 \sigma^2$}
   & & &        \\
\hline
4  \phv{11}{5}     
   & \onecoll{1}{$- 4.1759835935184954553812374$}
   & \onecoll{1}{$- 0.2689357165826222969605534$}
   & $\frac{1}{2} \sigma$
   & Yes
   & $2$        \\
   \phv{5}{5}        
   & \onecolr{1}{${}- 0.5378714331652445939211068 \sigma$}
   & \onecolr{1}{${}- 2.0879917967592477276906187 \sigma$}
   & & &        \\
\hline
5  \phv{11}{5}   
   & \onecoll{1}{$0.3647664814623851156843183$}
   & \onecoll{1}{$0.1044209421210358098319009$}
   & $\frac{1}{2} (-1  + \sigma - \sigma^2)$
   & Yes
   & $2$        \\
   \phv{5}{5}       
   & \onecol{1}{${}+ 0.9383748471668418510324386 \sigma$}
   & \onecol{1}{${}- 1.0427957892878776608643396 \sigma$}
   & & &        \\
   \phv{5}{5}        
   & \onecolr{1}{${}+ 1.5119832128712985863805590 \sigma^2$}
   & \onecolr{1}{${}- 0.4691874235834209255162193 \sigma^2$}
   & & &        \\
\hline
6  \phv{11}{5}   
   & \onecoll{1}{$0.3903032175535131898951365 (1 + \sigma^2)$}
   & \onecoll{1}{$0.6769888476508730716284981$}
   & $\frac{1}{2} (-1  - \sigma + \sigma^2)$
   & Yes
   & $2$        \\
   \phv{5}{5}        
   & \onecolr{1}{${}+ 2.1345841304087725230472693 \sigma$}
   & \onecolr{1}{${}- 1.0672920652043862615236346 (\sigma + \sigma^2)$}
   & & &        \\
\hline
7  \phv{11}{5}     
   & \onecoll{1}{$- 0.5430424606759486694736326 (1 + \sigma)$}
   & \onecoll{1}{$  1.0860849213518973389472652$}
   & $\frac{1}{5}(-3 + 2\sigma + 2\sigma^2 + 2\sigma^3 - 3\sigma^4)$
   & No
   & $1$        \\
   \phv{5}{5}       
   & \onecol{1}{${}+ 0.8649249218235797385747707 (\sigma^2 +
                                                       \sigma^4)$}
   & \onecol{1}{${}- 0.3218824611476310691011381 (\sigma +
                                                     \sigma^4$)}
   & & &          \\
   \phv{5}{5}        
   & \onecolr{1}{\qquad ${}- 0.6437649222952621382022762 \sigma^3$}
   & \onecolr{1}{\qquad ${}- 0.2211599995283176003724945 (\sigma^2 +
\sigma^3)$} 
   & & &        \\
\hline
8  \phv{11}{5}     
   & \onecoll{1}{$- 1.3494436538740630114284692$}
   & \onecoll{1}{$  0.4769621163349386255017257$}
   & \onecoll{1}{$\frac{1}{82} (44 + 29\sigma - 34\sigma^2 +
13\sigma^3 + 30\sigma^4 + 3\sigma^5$}
   & Yes
   & $2$        \\
   \phv{5}{5}       
   & \onecol{1}{${}- 1.7116028784482896128599650 \sigma$}
   & \onecol{1}{${}- 1.3960582877976217506952332 \sigma$}
   & \onecolr{1}{${}- 12\sigma^6 + 7\sigma^7 - 28\sigma^8 - 11\sigma^9)$}
   & &          \\
   \phv{5}{5}       
   & \onecol{1}{\ ${}- 0.3796444360048927748408130 \sigma^2$}
   & \onecol{1}{\ ${}- 0.6175505522167333606273630 \sigma^2$}
   & & &          \\
   \phv{5}{5}
   & \onecol{1}{\ \ ${}- 0.8340983749986581978191287 \sigma^3$}
   & \onecol{1}{\ \ ${}- 0.4964278267523899302068686 \sigma^3$}
   & & &        \\
   \phv{5}{5}       
   & \onecol{1}{\ \ \ ${}- 0.1806103940721755799158468 \sigma^4$}
   & \onecol{1}{\ \ \ ${}- 0.2291274044079765136578312 \sigma^4$}
   & & &        \\
   \phv{5}{5}
   & \onecol{1}{\ \ \ \ ${}+ 0.7063623627724616272858469 \sigma^5$}
   & \onecol{1}{\ \ \ \ ${}- 0.5269815315477649828362206 \sigma^5$}
   & & &        \\
   \phv{5}{5}
   & \onecol{1}{\ \ \ \ \ ${}- 1.2043707169764199179058313 \sigma^6$}
   & \onecol{1}{\ \ \ \ \ ${}+ 0.2599539512328980544134731 \sigma^6$}
   & & &        \\
   \phv{5}{5}       
   & \onecol{1}{\ \ \ \ \ \ ${}+ 1.0240929894938156671027431 \sigma^7$}
   & \onecol{1}{\ \ \ \ \ \ ${}- 0.1966713232003592665532588 \sigma^7$}
   & & &        \\
   \phv{5}{5}        
   & \onecol{1}{\ \ \ \ \ \ \ ${}- 2.4037576815736824020478056 \sigma^8$}
   & \onecol{1}{\ \ \ \ \ \ \ ${}- 0.0966340491363850966012588 \sigma^8$}
   & & &        \\
   \phv{5}{5}         
   & \onecolr{1}{${}+ 0.4024540247941539343479953 \sigma^9$}
   & \onecolr{1}{${}- 0.1427744719524809127778017 \sigma^9$}
   & & &        \\
\hline
\end{tabular}
}
\end{sidewaystable}

\begin{sidewaystable}
\centering
{\fontsize{8pt}{10.5pt}\selectfont
\begin{tabular}{|c|c|c|c|c|c|}
\hline
\onecolc{6}{\footnotesize \textrm{Table 5:} \emph{Verification of
Conjecture 2.2 - {\rm part (i), (iii) and (iv)}
\phv{11}{5}}} \\
\hline
\# \phv{11}{5}
   & \makebox[7cm][c]{$4\sqrt{d_k}\inv R(\gamma)$}
   & \makebox[7.5cm][c]{$\Phi_{\ff,\emptyset}(1)$}
   & \makebox[4cm][c]{$A$}
   & is $\ff = \fq^l$?
   & $d_\ff$        \\
\hline
9  \phv{11}{5}    
   & \onecoll{1}{$ 1.2006136993027519158356429$}
   & \onecoll{1}{$ 0.3555763402274556807627663$}
   & $\frac{1}{2} (-1  - \sigma + \sigma^2 - \sigma^3 - \sigma^4 +
\sigma^5 + \sigma^6)$
   & Yes
   & $2$        \\
   \phv{5}{5}       
   & \onecol{1}{${}+ 0.5601281150225709497028443 \sigma$}
   & \onecol{1}{${}- 0.8024573768066609557593419 \sigma$}
   & & &          \\
   \phv{5}{5}        
   & \onecol{1}{\ ${}+ 0.4468810365792052749965755 \sigma^2$}
   & \onecol{1}{\ ${}- 0.3366875967329683122045565 \sigma^2$}
   & & &        \\
   \phv{5}{5}
   & \onecol{1}{\ \ ${}+ 0.3336339581358396002903067 \sigma^3$}
   & \onecol{1}{\ \ ${}- 0.3981563224960909600763009 \sigma^3$}
   & & &        \\
   \phv{5}{5}
   & \onecol{1}{\ \ \ ${}- 0.3068516261443413658424917 \sigma^4$}
   & \onecol{1}{\ \ \ ${}- 0.2234405182896026374982877 \sigma^4$}
   & & &        \\
   \phv{5}{5}
   & \onecol{1}{\ \ \ \ ${}+ 0.9126508166528979185513609 \sigma^5$}
   & \onecol{1}{\ \ \ \ ${}- 0.0487247140831143149202745 \sigma^5$}
   & & &        \\
   \phv{5}{5}         
   & \onecolr{1}{${}- 0.0188887434944873685582098 \sigma^6$}
   & \onecolr{1}{${}- 0.1101934398462369627920190 \sigma^6$}
   & & &        \\
\hline
10 \phv{11}{5}     
   & \onecoll{1}{$ 0.7234352393016752990818922 (1 + \sigma)$}
   & \onecoll{1}{$ 0.3662434964105941418610148$}
   & $\frac{1}{2} (-1  + \sigma - \sigma^2)$
   & Yes
   & $2$        \\
   \phv{5}{5}         
   & \onecolr{1}{${}+ 2.1793574714245388818858141 \sigma^2$}
   & \onecolr{1}{${}- 1.0896787357122694409429070 (\sigma + \sigma^2)$}
   & & &        \\
\hline
11 \phv{11}{5}     
   & \onecoll{1}{$ 0.7895116790170794653416756$}
   & \onecoll{1}{$ 2.6211847478167423772732258$}
   & $\frac{1}{6}(-1 - \sigma + 5\sigma^2 - \sigma^3 - \sigma^4 - \sigma^5)$
   & No
   & $1$        \\
   \phv{5}{5}        
   & \onecol{1}{${}- 1.1652625156938787944693039 \sigma$}
   & \onecol{1}{${}+ 0.3526417361775253842417098 \sigma$}
   & & &          \\
   \phv{5}{5}
   & \onecol{1}{\ ${}- 1.8085639683003889670456318 \sigma^2$}
   & \onecol{1}{\ ${}+ 0.7895116790170794653416756 \sigma^2$}
   & & &        \\
   \phv{5}{5}        
   & \onecol{1}{\ \ ${}- 0.7895116790170794653416756 \sigma^3$}
   & \onecol{1}{\ \ ${}- 1.1652625156938787944693039 \sigma^3$}
   & & &        \\
   \phv{5}{5}
   & \onecol{1}{\ \ \ ${}+ 2.6211847478167423772732258 \sigma^4$}
   & \onecol{1}{\ \ \ ${}- 1.8085639683003889670456318 \sigma^4$}
   & & &        \\
   \phv{5}{5}         
   & \onecolr{1}{${}+ 0.3526417361775253842417098 \sigma^5$}
   & \onecolr{1}{${}- 0.7895116790170794653416756 \sigma^5$}
   & & &        \\
\hline
12 \phv{11}{5}       
   & \onecoll{1}{$-  0.2877586687247090106420884$}
   & \onecoll{1}{${} 0.1438793343623545053210442$}
   & $-\frac{1}{2}$
   & Yes
   & $2$        \\
   \phv{5}{5}         
   & \onecolr{1}{$  + 3.9680836522391984256974575 \sigma$}
   & \onecolr{1}{${}- 1.9840418261195992128487287 \sigma$}
   & & &        \\
\hline
13 \phv{11}{5}     
   & \onecoll{1}{$ 2.3349046592276594288814979$}
   & \onecoll{1}{$ 1.3358165556295184849328475$}
   & $\frac{1}{2} (1  - \sigma - \sigma^2)$
   & Yes
   & $2$        \\
   \phv{5}{5}         
   & \onecolr{1}{${}- 1.1674523296138297144407489 (\sigma + \sigma^2)$}
   & \onecolr{1}{${}- 0.1683642260156887704920986 (\sigma + \sigma^2)$}
   & & &        \\
\hline
14 \phv{11}{5}      
   & \onecoll{1}{$  0.4165560795603681099566393$}
   & \onecoll{1}{$- 0.4165560795603681099566393$}
   & $\frac{1}{6}(-5 + \sigma^3)$
   & No
   & $1$        \\
   \phv{5}{5}        
   & \onecol{1}{${}+ 0.7153308793910213046570109 \sigma$}
   & \onecol{1}{${}- 0.7153308793910213046570109 \sigma$}
   & & &        \\
   \phv{5}{5}
   & \onecol{1}{\ ${}- 1.2682208798620716367535322 \sigma^2$}
   & \onecol{1}{\ ${}+ 1.2682208798620716367535322 \sigma^2$}
   & & &        \\
   \phv{5}{5}
   & \onecol{1}{\ \ ${}- 0.4165560795603681099566393 \sigma^3$}
   & \onecol{1}{\ \ ${}+ 0.4165560795603681099566393 \sigma^3$}
   & & &        \\
   \phv{5}{5}        
   & \onecol{1}{\ \ \ ${}- 0.7153308793910213046570109 \sigma^4$}
   & \onecol{1}{\ \ \ ${}+ 0.7153308793910213046570109 \sigma^4$}
   & & &        \\
   \phv{5}{5}         
   & \onecolr{1}{${}+ 1.2682208798620716367535322 \sigma^5$}
   & \onecolr{1}{${}- 1.2682208798620716367535322 \sigma^5$}
   & & &        \\
\hline
15 \phv{11}{5}      
   & \onecoll{1}{$- 1.2218884525300797394860547 (1 + \sigma_1^2)$}
   & \onecoll{1}{$  2.9561452818637877653536663$}
   & \onecoll{1}{$\frac{1}{18} (1 + \sigma_1 + \sigma_1^2 - 11\sigma_2
+ 7 \sigma_1\sigma_2$} 
   & Yes
   & $2$        \\   
   & \onecol{1}{${}+ 3.4685136586674160517352232 \sigma_1$}
   & \onecol{1}{${}- 1.7342568293337080258676116
                     (\sigma_1 + \sigma_1^2)$}
   & \onecolr{1}{${}+ 7\sigma_1^2\sigma_2 + \sigma_2^2 +
\sigma_1\sigma_2^2 + \sigma_1^2\sigma_2^2)$} 
   & &          \\
   \phv{5}{5}        
   & \onecol{1}{\ ${}- 0.6973303591485227158417690
                     (\sigma_2 + \sigma_1^2\sigma_2^2)$}
   & \onecol{1}{\ ${}- 1.0705305164860507528782365
                     (\sigma_2 + \sigma_2^2)$}
   & & &          \\
   \phv{5}{5}
   & \onecol{1}{$\ \ {}+ 2.2802292524382017551015625
                     (\sigma_1^2\sigma_2 + \sigma_2^2)$}
   & \onecol{1}{$\ \ {}+ 1.7678608756345734687200056
                     (\sigma_1\sigma_2 + \sigma_1^2\sigma_2^2)$}
   & & &          \\
   \phv{5}{5}            
   & \onecolr{1}{~~~${}- 0.5581621396824224664966796
                         (\sigma_1\sigma_2 + \sigma_1\sigma_2^2)$}
   & \onecolr{1}{~~~${}- 1.2096987359521510022233259
                         (\sigma_1^2\sigma_2 + \sigma_1\sigma_2^2)$}
   & & &        \\
\hline
\end{tabular}
}
\end{sidewaystable}

\begin{sidewaystable}
\centering
{\fontsize{9pt}{12pt}\selectfont
\begin{tabular}{|c|c|c|c|c|}
\hline
\onecolc{5}{\footnotesize \textrm{Table 6:} \emph{Verification of
Conjecture 2.2 - {\rm part (ii) and (iv)}
\phv{11}{5}}} \\
\hline
\# \phv{11}{5}
   & \makebox[2cm][c]{primes}
   & \makebox[14cm][c]{$C_{p,\ff} R_p (\eta_\ff) = \Phi_{\ff, T_p, p}(1)$}
   & \makebox[1.5cm][c]{$d_{\ff, \sigma - 1}$}
   & \makebox[3cm][c]{index of $\ZG\eta_\ff$}
   \\
\hline
1  \phv{11}{5}
   & $3$, $7$, $11$
   & \onecoll{1}{$0.202021222001202022011122201001212121201_{3}$}
   & $1$
   & $2^2$
   \\
   &
   & \onecolr{1}{${}+ 0.002112222212110120202010210110011000011_{3}
     (\sigma + \sigma^2)$}
   & &          \\
   \phv{11}{5}
   &
   & $0.232034003422155306164163_{7} +
      0.624214462041162660106331_{7} (\sigma + \sigma^2)$
   & &          \\
   \phv{11}{5}
   &
   & $0.859AA8491A4592272_{11} + 0.593A1A1A496337044_{11}
      (\sigma + \sigma^2)$
   & &          \\
\hline
2  \phv{11}{5}
   & $19$
   & $0.37A267_{19} + 0.IFBDF1_{19} \sigma + 0.7C1858_{19} \sigma^2$
   & $1$
   & $2^2$
   \\
\hline
3  \phv{11}{5}
   & $3$, $11$
   & \onecoll{1}{$0.022000201100100201021001020122020221201_{3}$}
   & $1$
   & $2^2$ 
   \\
   \phv{5}{5}
   &
   & \onecolc{1}{${}+ 0.011221101110201021102020011101021020122_{3}
   (\sigma + \sigma^3)$}
   & & \\
   \phv{5}{5}
   &
   & \onecolr{1}{${}+ 0.000120101120012102212010002112212112201_{3}
   \sigma^2$} 
   & & \\
   \phv{11}{5}
   &
   & $0.5281109A901147AA7_{11} + 0.065022A7402839018_{11} (\sigma +
\sigma^3) + 0.692A2405AA2430228_{11} \sigma^2$
   & &          \\
\hline
4  \phv{11}{5}
   & $11$
   & $0.9627359501683452_{11} + 0.00637222A6760375_{11} \sigma$
   & $1$
   & $2$
   \\
\hline
5  \phv{11}{5}
   & $5$, $13$
   & $0.44203011304041240231402_{5}
      + 0.41443324422220142233001_{5} \sigma
      + 0.40101242140401414400031_{5} \sigma^2$
   & $1$
   & $2^2$
   \\
   \phv{11}{5}
   &
   & $0.5811AA04_{13} + 0.C408A99C_{13} \sigma + 0.201B3AB1_{13} \sigma^2$
   & &          \\
\hline
6  \phv{11}{5}
   & $3$, $17$
   & \onecoll{1}{$0.102212201210222202002020222122212022010_{3}$}
   & $1$
   & $2^2$
   \\
   \phv{5}{5}
   &
   & \onecolr{1}{$+ 0.011001110112222111101011110202122110222_{3}
      (\sigma + \sigma^2)$}
   & &          \\
   \phv{11}{5}
   &
   & $0.F44A49F_{17} + 0.8218BE5_{17} (\sigma + \sigma^2)$
   & &          \\
\hline
7  \phv{11}{5}
   & $5$, $7$, $11$
   & \onecoll{1}{$0.3233210340311430413102022_{5}
      + 0.0311412013321131241321400_{5} (\sigma + \sigma^4)$}
   & $1$
   & $1$
   \\
   \phv{5}{5}
   &
   & \onecolr{1}{${}+ 0.1314122043432120343022033_{5} (\sigma^2 + \sigma^3)$}
   & &          \\
   \phv{11}{5}
   &
   & \onecoll{1}{$0.3026023225325560324160532_{7}
      + 0.0165614661060504553104546_{7} (\sigma + \sigma^4)$}
   & &          \\
   \phv{5}{5}
   &
   & \onecolr{1}{${}+ 0.2231301655466522334432206_{7} (\sigma^2 + \sigma^3)$}
   & &          \\
   \phv{11}{5}
   &
   & $0.329784702951_{11} + 0.889019601695_{11} (\sigma + \sigma^4)
      + 0.7670050A8599_{11} (\sigma^2 + \sigma^3)$
   & &          \\
\hline
8  \phv{11}{5}
   & $11$, $41$
   & \onecoll{1}{$0.8806785A3211A8230_{11}
      + 0.06600615342772379_{11} \sigma
      + 0.60514793493396165_{11} \sigma^2$}
   & $1$
   & $2^9$
   \\
   \phv{5}{5}
   &
   & \onecolc{1}{${}+ 0.05A986A29673A92A3_{11} \sigma^3
      + 0.3496587A010380A79_{11} \sigma^4$}
   & &          \\
   \phv{5}{5}
   &
   & \onecolc{1}{\ ${}+ 0.49400712781309175_{11} \sigma^5
      + 0.99099496585661127_{11} \sigma^6$}
   & &          \\
   \phv{5}{5}
   &
   & \onecolr{1}{${}+ 0.6637455738542124A_{11} \sigma^7
      + 0.5619A049053576157_{11} \sigma^8
      + 0.60603562051742618_{11} \sigma^9$}
   & &          \\
   \phv{11}{5}
   &
   & \onecoll{1}{$0.SI_{41} + 0.5(37)_{41} \sigma + 0.BL_{41} \sigma^2
      + 0.1Q_{41} \sigma^3 + 0.3W_{41} \sigma^4 + 0.WL_{41} \sigma^5$}
   & &          \\
   \phv{5}{5}
   &
   & \onecolr{1}{${}+ 0.Y(36)_{41} \sigma^6 + 0.VN_{41} \sigma^7
                  + 0.SA_{41} \sigma^8 + 0.D(37)_{41} \sigma^9$}
   & &          \\
\hline
\end{tabular}
}
\end{sidewaystable}

\begin{sidewaystable}
\centering
{\fontsize{9pt}{12pt}\selectfont
\begin{tabular}{|c|c|c|c|c|}
\hline
\onecolc{5}{\footnotesize \textrm{Table 7:} \emph{Verification of
Conjecture 2.2 - {\rm part (ii) and (iv)}
\phv{11}{5}}}
\\
\hline
\# \phv{11}{5}
   & \makebox[2cm][c]{primes}
   & \makebox[14cm][c]{$C_{p,\ff} R_p (\eta_\ff) = \Phi_{\ff, T_p, p}(1)$}
   & \makebox[2cm][c]{$d_{\ff, \sigma - 1}$}
   & \makebox[3cm][c]{index of $\ZG\eta_\ff$}
   \\
\hline
9  \phv{11}{5}
   & $3$, $11$, $17$
   & \onecoll{1}{$0.222022110111002222112202221101101011211_{3}
      + 0.002001201222221221212201112021202201020_{3} \sigma$}
   & $1$
   & $2^6$
   \\
   \phv{5}{5}
   &
   & \onecolc{1}{${}+ 0.002001220120020121010102102021122102221_{3} \sigma^2$}
   & &          \\
   &
   & \onecolc{1}{~${}+ 0.001211110201122210102210100110102201221_{3} \sigma^3$}
   & &          \\
   \phv{5}{5}
   &
   & \onecolc{1}{~~${}+ 0.112010220020002220212200202222100011101_{3} \sigma^4$}
   & &          \\
   \phv{5}{5}
   &
   & \onecolc{1}{~~~${}+ 0.220202201100201201000000011112201011010_{3} \sigma^5$}
   & &          \\
   \phv{5}{5}
   &
   & \onecolr{1}{${}+ 0.222022120210010100121111202101210110010_{3} \sigma^6$}
   & &          \\
   \phv{11}{5}
   &
   & \onecoll{1}{$0.78817332A270_{11} + 0.581648097113_{11} \sigma
     + 0.36370A254025_{11} \sigma^2 + 0.234186663552_{11} \sigma^3$}
   & &          \\
   \phv{5}{5}
   &
   & \onecolr{1}{${}+ 0.68A306753247_{11} \sigma^4
                   + 0.A26635843A21_{11} \sigma^5
                   + 0.9A6002162469_{11} \sigma^6$}
   & &          \\
   \phv{11}{5}
   &
   & \onecoll{1}{$0.092998_{17} + 0.411756_{17} \sigma
     + 0.BC80F6_{17} \sigma^2 + 0.DA572A_{17} \sigma^3$}
   & &          \\
   \phv{5}{5}
   &
   & \onecolr{1}{${}+ 0.25AGFF_{17} \sigma^4 + 0.8GE8C4_{17} \sigma^5
     + 0.AEBFG7_{17} \sigma^6$}
   & &          \\
\hline
10 \phv{11}{5}
   & $3$, $5$, $7$, $11$
   & \onecoll{1}{$0.122120101100101201201102002111221100102_{3}$}
   & $1$
   & $2^2$
   \\
   &
   & \onecolr{1}{$+ 0.010102011020202112122100011210200211221_{3}
     (\sigma + \sigma^2)$}
   & &          \\
   \phv{11}{5}
   &
   & $0.40300320211333232201403121_{5}
      + 0.24142112410003431344112113_{5} (\sigma + \sigma^2)$
   & &          \\
   \phv{11}{5}
   &
   & $0.613562234546646014324320_{7}
      + 0.644312242351243462416504_{7} (\sigma + \sigma^2)$
   & &          \\
   \phv{11}{5}
   &
   & $0.89816292686A4601_{11} + 0.20928A9982220855_{11} (\sigma + \sigma^2)$
   & &          \\
\hline
11 \phv{11}{5}
   & $11$
   & \onecoll{1}{$0.4279553895127000_{11}
      + 0.3502739301069659_{11} \sigma + 0.4778387406583A34_{11} \sigma^2$}
   & $1$
   & $1$
   \\
   \phv{5}{5}
   &
   & \onecolr{1}{${}+ 0.273038669A268856_{11} \sigma^3
      + 0.27A95423246775A5_{11} \sigma^4 + 0.73327236A4527076_{11} \sigma^5$}
   & &          \\
\hline
12 \phv{11}{5}
   & $11$, $41$
   & $0.23A9227A0541A025_{11} + 0.25A583269A5537A1_{11} \sigma$
   & $1$
   & $2$
   \\
   \phv{11}{5}
   &
   & $0.6SP_{41} +  0.5KN_{41} \sigma$
   & &          \\
\hline
13 \phv{11}{5}
   & $7$, $11$
   & $0.5605361401660032563563_{7}
      + 0.1610430464366153334222_{7} (\sigma + \sigma^2)$
   & $1$
   & $2^2$
   \\
   \phv{11}{5}
   &
   & $0.2398796701A1346A_{11}
      + 0.288A2880995A2406_{11} (\sigma + \sigma^2)$
   & &          \\
\hline
14 \phv{11}{5}
   & $11$
   & \onecoll{1}{$0.290A17A7368678838_{11}
                  + 0.8A649AA3287793065_{11} \sigma
                  + 0.7663619146A559404_{11} \sigma^2$}
   & $1$
   & $1$
   \\
   \phv{5}{5}
   &
   & \onecolr{1}{$+ 0.91A09303742432272_{11} \sigma^3
                  + 0.30461007823317A45_{11} \sigma^4
                  + 0.444749196405516A6_{11} \sigma^5$}
   & &          \\
\hline
15 \phv{11}{5}
   & $13$, $19$
   & \onecoll{1}{$0.37861615758_{13}
      + 0.B7962270777_{13} (\sigma_1 + \sigma_1^2)
      + 0.6853599605C_{13} (\sigma_2 + \sigma_2^2)$}
   & $d_{\ff, \sigma_1 - 1} = 1$
   & $2^6$
   \\
   &
   & \onecolr{1}{${}+ 0.09A939C0279_{13}
     (\sigma_1\sigma_2 + \sigma_1^2\sigma_2^2)
     + 0.65B6A46B581_{13} (\sigma_1^2\sigma_2 + \sigma_1\sigma_2^2)$}
   & $d_{\ff, \sigma_2 - 1} = 1$
   &            \\
   \phv{13}{5}
   &
   & \onecoll{1}{$0.BC66483F_{19} + 0.G0B815H6_{19} (\sigma_1 + \sigma_1^2)
     + 0.7BB0F5H3_{19} (\sigma_2 + \sigma_2^2)$}
   & &          \\
   \phv{5}{5}
   &
   & \onecolr{1}{${}+ 0.66709GBE_{19} (\sigma_1\sigma_2
     + \sigma_1^2\sigma_2^2) + 0.BF932F8A_{19}
     (\sigma_1^2\sigma_2 + \sigma_1\sigma_2^2)$}
   & &          \\
\hline
\end{tabular}
}
\end{sidewaystable}

\end{document}